\documentclass[12pt]{amsart}
\usepackage[matrix,arrow,curve]{xy}
\usepackage{amssymb}
\usepackage{a4}

\newtheorem{dummy}{anything}[section] 
\newtheorem{Theorem}[dummy]{Theorem}
\newtheorem*{thma}{Theorem A}
\newtheorem*{thmb}{Theorem B}
\newtheorem{Lemma}[dummy]{Lemma} 
\newtheorem{Proposition}[dummy]{Proposition}
 \newtheorem*{propc}{Proposition C}
\newtheorem*{propd}{Proposition D}
\newtheorem*{prope}{Proposition E} 
 \newtheorem{Corollary}[dummy]{Corollary} 
\theoremstyle{definition}               

\newtheorem{Example}[dummy]{Example}
\newtheorem{Remark}[dummy]{Remark}
\newtheorem{ccote}[dummy]{}

\newcommand{\preu}{\begin{proof}}

\newcommand{\fl}[1]{\buildrel{#1}\over{\longrightarrow}}
\newcommand{\hfl}[2]{\smash{\mathop{\hbox to 1 truecm{\kern
3pt\rightarrowfill\kern 3pt}}%
\limits^{\scriptstyle#1}_{\scriptstyle#2}}}
\newcommand{\cqfd}{\end{proof}}
\newcommand{\equationref}[1]{(\ref{#1})}
\newcommand{\bbr}{{\mathbf R}}

\newcommand{\bbz}{{\mathbf Z}}

\newcommand{\bbn}{{\mathbf N}}
\newcommand{\bbh}{{\mathbf H}}

\newcommand{\cale}{{\mathcal E}}
\newcommand{\calf}{{\mathcal F}}

\newcommand{\calh}{{\mathcal H}}

\newcommand{\calr}{{\mathcal R}}
\newcommand{\cals}{{\mathcal S}}
\newcommand{\calt}{{\mathcal T}}

\newcommand{\pcirc}{\kern .7pt {\scriptstyle \circ} \kern 1pt}
\newcommand{\iso}{\cong}
\newcommand{\mun}{{-1}}
\newcommand{\vv}{{\,|\,}}
\newcommand{\res}{{\rm Res\,}}
\newcommand{\aut}{{\rm Aut}}
\newcommand{\smallpm}{{\scriptscriptstyle \pm}}
\newcommand{\smallplus}{{\scriptscriptstyle +}}
\newcommand{\smallminus}{{\scriptscriptstyle -}}

\newcommand{\eqncount}{\setcounter{equation}{\value{dummy}}%
\addtocounter{dummy}{1}}

\begin{document}
\title[Equivariant principal bundles]
{Equivariant principal bundles over spheres\\ and cohomogeneity one
manifolds}
\author{Ian HAMBLETON} 
\thanks{\noindent The authors thank the Swiss National Fund 
for Scientific Research, the Universit\'e de G\'en\`eve, and the Max Planck Institut f\"ur Mathematik in Bonn for
hospitality and support.}
\address{Department of Mathematics \& Statistics
\newline\indent
McMaster University
\newline\indent
Hamilton, ON  L8S 4K1, Canada}
\email{ian{@}math.mcmaster.ca}
\author{Jean-Claude HAUSMANN}
\address{Section de Math\'ematiques
\newline\indent
Universit\'e de Gen\`eve,
 B.P. 240
\newline\indent
CH-1211 Gen\`eve 24,
Switzerland}
\email{
hausmann{@}math.unige.ch}
\date{Jan. 20, 2001}  
\begin{abstract}\noindent 
We classify $SO(n)$-equivariant principal bundles 
over $S^n$ in terms of their isotropy 
representations over the north and south poles.
This is an example of a general result classifying equivariant
$(\Pi, G)$-bundles over cohomogeneity one manifolds.
\end{abstract}
\maketitle

\section{Introduction}\label{intro1}

Let $\Pi$ and $G$ be Lie groups.
A principal $(\Pi,G)$-bundle is a locally trivial,
principal $G$-bundle $p\colon E\to X$ such that $E$ and $X$ 
are left $\Pi$-spaces. The projection map $p$ is $\Pi$-equivariant 
and $\gamma(e\cdot g)= (\gamma e)\cdot g$, where $\gamma \in \Pi$ and
$g\in G$ acts on $e\in E$  by the principal action.
Equivariant principal bundles, and their natural generalizations, were 
studied by T. E. Stewart \cite{S},
T. tom Dieck \cite{D1}, \cite[I\,(8.7)]{D2}, R. Lashof
\cite{L2}, \cite{L3} together with P. May \cite{LM} and G. Segal \cite{LMS}.

These authors study equivariant principal bundles by homotopy
theoretic methods.
There exists a classifying space $B(\Pi, G)$  for principal  
$(\Pi,G)$-bundles \cite{D1}, 
so the classification of equivariant bundles in
particular cases can be approached by studying the $\Pi$-equivariant
homotopy type of $B(\Pi, G)$. If the structural group $G$ of the bundle
is {\it abelian}, then the main result of \cite{LMS} states that equivariant
bundles over a $\Pi$-space $X$ are classified by the ordinary
homotopy classes of maps $[X \times_{\Pi} E\Pi, BG]$.
In practice, this program leads to an obstruction
theory rather than a classification. 
See, however, the results of Lashof in the special cases where
$\Pi$ acts transitively \cite{L1} or semi-freely \cite{L3} on
the base space $X$.
 
Another approach to equivariant principal bundles uses the
``local'' invariants arising from isotropy representations at
singular points of $(X,\Pi)$,
together with equivariant gauge theory \cite{BH}, \cite{HL}, 
\cite{HLM}, \cite{Ha}. 
By an isotropy representation at a $\Pi$-fixed point $x_0 \in X$ 
we mean the homomorphism $\alpha_{x_0}\colon\Pi\to G$
defined by the formula 
$$\gamma\cdot e_0 = e_0 \cdot \alpha(\gamma)$$
where $e_0 \in p^{-1}(x_0)$. The homomorphism $\alpha$ is independent
of the choice of $e_0$ up to conjugation in $G$.
The relationship between the local invariants and the homotopy
classification (in the form of a
Localization Theorem ?) deserves further study.

In this paper, we use the second approach for
$\Pi = SO(n)$ acting in the standard way on $X = S^n$.
In this concrete situation, we obtain a complete classification by
relatively elementary geometric methods. It turns out that
the local isotropy representations at the north and south poles of
$S^n$ explicitly determine the classification of $(SO(n), G)$
principal bundles over $S^n$ (for short $(n,G)$-bundles). 
One surprising consequence is
that the set $\cale (n, G)$ of $(n,G)$-bundles is finite for $n\geq 3$.
In contrast, the set of (non-equivariant) principal 
$G$-bundles over $S^n$ is often infinite.
 A detailed statement of these
results is given in the next section and their proofs, essentially
self-contained, are explained  in Sections 
\ref{preli} to \ref{S:son+1}.
Several examples are given in Section \ref{exap}.
In Section \ref{genset} we show how these results fit into the 
more general setting of equivariant $(\Pi,G)$-bundles
over certain $\Pi$-manifolds studied by K.~J\"anich \cite{Ja}
and E.~Straume \cite{Straume}.
In particular, we obtain a classification of $(\Pi,G)$-bundles
over cohomogeneity 1 manifolds.

The authors would like to thank P. de la Harpe for helpful discussions.

\section{Statement of results}\label{intro}

Let $S^n$ be the $n$-dimensional sphere of radius 1 in $\bbr^{n+1}$.
We consider the action on $S^n$ of the group $SO(n)$, 
by orthogonal transformations
fixing the poles $(0,\dots ,0,\pm 1)$.

Let $G$ be a Lie group. 
We denote by $\calr(n,G)$ the set 
of smooth homomorphisms from $SO(n)$ to $G$
modulo the conjugations by elements of $G$. Unless specified,
all maps between manifolds are smooth of class $C^\infty$.

By a $G$-principal bundle $\eta$ over $S^n$, 
we mean, as usual, a smooth map $p\colon  E\to S^n$ from 
a manifold $E=E(\eta)$ and a free right action $E\times G\to E$
so that $p(z\cdot g) = p(z)$ with the standard local triviality condition.
The isomorphism classes of $G$-bundles over $S^n$ are
in bijection with $\pi_{i-1}(G)/\pi_0(G)$, the quotient of the homotopy group
$\pi_{n-1}(G)$ (based on the neutral element $e$ of $G$) by the action of
$\pi_0(G)$
induced by the conjugation of $G$ on itself.
The bijection associates to a bundle $\eta$ the class 
$C(\eta):=[\partial({\rm id\,}_{S^n})]\in \pi_{n-1}(G)/\pi_0(G)$, where
$\partial \colon  \pi_n(S^n)\to\pi_{n-1}(G)$ is the boundary operator in the
homotopy exact
sequence of $\eta$ \cite[Th.\, 18.5]{St}. 

A {\it $SO(n)$-equivariant principal $G$-bundle} 
$\xi$ over $S^n$ (or an {\it $(n,G)$-bundle} for short) 
is a $G$-principal bundle $\xi^\flat$ over $S^n$
together with a left action $SO(n)\times E(\xi)\to E(\xi)$ 
commuting with the free right action of $G$ and  
such that the projection to $S^n$ is $SO(n)$-equivariant
(we write $E(\xi)$ for $E(\xi^\flat$)). 
Two $(n,G)$-bundles $\xi_1,\xi_2$ are {\it isomorphic} if there exists 
a  diffeomorphism $h: E(\xi_1)\to E(\xi_2)$ 
which is both $SO(n)$ and $G$-equivariant and which commutes
with the projections to $S^n$.
We will compute the set $\cale (n,G)$
of isomorphism classes of $(n,G)$-bundles.

Let $\xi$ be a $(n,G)$-bundle. 
Choose points $a,b\in E(\xi)$ such that 
$p(a)=(0,\dots ,-1)$ and $p(b)=(0,\dots ,1)$.
Let $\alpha,\beta$ be the maps from $SO(n)$ to $G$ 
determined by the formulae
$A\cdot a = a \cdot \alpha(A)$ and $A\cdot b = b \cdot \beta(A)$.
We shall prove in Lemma \ref{homos} that 
$\alpha$ and $\beta$ are smooth homomorphisms and that their
class in $\calr(n,G)$ depend only on $[\xi]\in\cale (n,G)$.
We call $\alpha$ and
$\beta$ the {\it isotropy representations} (associated to
$a$ and $b$).
This defines a map $J\colon \cale (n,G)\to\calr(n,G)\times\calr(n,G)$ 
by $J(\xi):=([\alpha],[\beta])$. 
When $n=2$ and $G$ is connected, 
$J(\xi)$ is a complete invariant which, 
in particular, determines
the (non-equivariant) isomorphism class of $\xi^\flat$. 
More precisely, let 
$\psi\colon \calr(2,G)\times\calr(2,G)\to \pi_1(G)$ 
be the map determined by 
$\psi(\alpha,\beta)(z):=[\alpha(z)\beta(z)^\mun]$
(one uses that $SO(2)\approx S^1$ and that 
$\psi$ is well defined if $G$ is connected).

\begin{thma}
Suppose that $G$ is connected Lie group. Then,
\begin{enumerate}
\renewcommand{\labelenumi}{(\roman{enumi})}
\item the map 
$J\colon  \cale(2,G)\to \calr(2,G)\times\calr(2,G)$ is a bijection.
\item if $J(\xi)=([\alpha],[\beta])$, then 
$\psi(\alpha,\beta)=C(\xi^\flat)$.
\end{enumerate}
\end{thma}


We shall now generalize Theorem A for $n\geq 2$ or $G$ any Lie group.
In general, $J$ is then neither injective nor surjective and $C(\xi^\flat)$
is not determined by $J(\xi)$. 
Consider $SO(n-1)$ as the subgroup of 
$SO(n)$ fixing the last 
coordinate. The restriction 
$[\mu]\mapsto [\mu\, {\scriptstyle |}_{SO(n-1)}]$ gives a map
$\res\colon\calr(n,G)\to\calr(n-1,G)$. Denote by 
$\calr(n,G)\times_{(n-1)}\calr(n,G)$ the set of 
$([\alpha],[\beta])\in\calr(n,G)\times\calr(n,G)$
such that $\res[\alpha]=\res[\beta]$.
If $\varphi\colon H\to G$ is a group homomorphism,
we denote by $Z_\varphi\subset G$ the centralizer
of $\varphi(H)$. 

\begin{thmb}
Let $G$ be any Lie group $G$. Then
\begin{enumerate}
\renewcommand{\labelenumi}{(\roman{enumi})}
\item the image of $J$ is $\calr(n,G)\times_{(n-1)}\calr(n,G)$.
\item let $\alpha,\beta\colon SO(n-1)\to G$ be two smooth homomorphisms
such that $\alpha\, {\scriptstyle |}_{SO(n-1)}=
\beta\, {\scriptstyle |}_{SO(n-1)}=\colon \gamma$. Then
$J^\mun([\alpha],[\beta])$ is in bijection with 
with the set of double cosets
$\pi_0(Z_\alpha)\backslash\pi_0(Z_\gamma)/\pi_0(Z_\beta)$.
\end{enumerate} 
\end{thmb}


\begin{Remark} The compatibility statement in Part (i) of Theorem~B 
was also observed by K.~Grove and W.~Ziller \cite[Prop.\,1.6]{GZ}. 
In \S\,\ref{genset}, Theorem B is  extended to a more general setting,
to include equivariant principal bundles over ``special'' $\Pi$-manifolds
in the sense of J\"anich \cite{Ja}. In particular this provides a
classification of the equivariant bundles considered by Grove and Ziller. 
\end{Remark}

Since $SO(1)$ is trivial, 
Theorem $B$ reduces to Part (i) of Theorem A when $n=2$.
To determine $C(\xi^\flat)$ as in Part (ii) of Theorem A, 
we must
choose particular representatives of $[\alpha]$ and $[\beta]$
(in general, $J(\xi)$ does not determine $\xi^\flat$: 
see examples \ref{2O2} and \ref{2kSOn}).
An {\it isotropic lifting} for $\xi$ is
a smooth curve $\tilde c\colon [-1,1]\to E(\xi)$ lifting
the meridian arc $c(t)=(0,\dots,\cos(\pi t/2),\sin(\pi t/2))$ and such that
$B\cdot\tilde c(t)=\tilde c(t)\alpha(B)$ for all $B\in SO(n-1)$.
Isotropic lifting always exist
(see Lemma \ref{tildec}). Choosing $a:=\tilde c(-1)$ and $b:=\tilde c(-1)$ leads
to isotropy representations $\alpha,\beta \colon  SO(n)\to G$
such that $\alpha\, {\scriptstyle |}_{SO(n-1)}=
\beta\, {\scriptstyle |}_{SO(n-1)}$. 
The map $\psi(\alpha,\beta)\colon SO(n)\to G$ 
constructed as in Theorem A then satisfies 
$\psi(\alpha,\beta)(AB)=\psi(\alpha,\beta)(A)$
when $B\in SO(n-1)$. It thus induces a map 
$$\bar\psi(\alpha,\beta)\colon  S^{n-1}\iso SO(n)/SO(n-1) \to G.$$
Note that $\bar\psi$ is well defined since
$\alpha$ and $\beta$ are actual homomorphisms and not conjugacy classes.

\begin{propc}
Let $\xi$ be a $(n,G)$-bundle. Let 
$\alpha,\beta \colon  SO(n)\to G$ be the isotropy representation
associated to the end points of an isotropic lifting.
Then, $[\bar\psi(\alpha,\beta)]=C(\xi^\flat)$
in $\pi_{n-1}(G)/\pi_0(G)$.
\end{propc}

We shall prove two consequences of Theorem B and Proposition C which emphasize
the constrast between the cases $n=2$ and $n\geq 3$.

\begin{propd} 
Let $\eta$ be a principal $G$-bundle over $S^2$
with $G$ a non-trivial Lie group. Then there exist infinitely
many $\xi\in\cale(2,G)$ such that $\xi^\flat\iso\eta$.
\end{propd}

\begin{prope} 
For $G$ a compact Lie group, the set $\cale (n,G)$ is finite when $n\geq 3$.
\end{prope}

These results are proved in 
\S\ \ref{S:pmr}, while the former sections are devoted to
preliminary material. In Section \ref{S:son+1}, we determine which
$(n,G)$-bundles come from an $SO(n+1)$-equivariant bundles.
Examples are  given in \S\,\ref{exap}.

\section{Preliminary constructions}\label{preli}

\begin{ccote}\label{Hwelldef} {\it $J$ is well defined.\/}
This follows from the following lemma.

\begin{Lemma}\label{homos}
Let $\xi$ be a $(n,G)$-bundle. 
Let $a,b\in E(\xi)$ such that 
$p(a)=(0,\dots ,-1)$ and $p(b)=(0,\dots ,1)$.
Let $\beta,\alpha$ be the maps from $SO(n)$ to $G$ 
determined by the formulae
$A\cdot a = a \cdot \alpha(A)$ and $A\cdot b = b \cdot \beta(A)$.
Then $\alpha$ and $\beta$ are smooth homomorphisms 
and their class in $\calr(n,G)$ depends only on $[\xi]\in\cale (n,G)$.
\end{Lemma}

\preu Let $A,B\in SO(n)$.
One has 
$$\begin{array}{rcl}
a\cdot \alpha(BA) &=& (BA)\cdot a = B\cdot (A\cdot a) = B\cdot (a\cdot
\alpha(A)) 
=\\[2pt] &=&
(B\cdot a) \cdot \alpha(A) = a\cdot (\alpha(B)\alpha(A)).
\end{array}$$
Therefore, $\alpha$ and similarly, $\beta$, are homomorphism. 
They are smooth because
the action of $SO(n)$ is smooth. 
If $a'$ is another choice instead of $a$,
there exists $g\in G$ such that $a'=a\cdot g$ and one has
\eqncount
\begin{equation}\label{eqconju}
a\cdot (g\alpha'(A)) = a'\cdot \alpha'(A) = A\cdot a' = A\cdot a\cdot g
= a\cdot (\alpha(A) g),
\end{equation}
whence $\alpha'(A)=g^\mun \alpha(A) g$.
This proves that the class of $(\alpha,\beta)$ in $\calr(n,G)$ 
does not depend on the choice of $a$ and $b$. Now, if 
$h\colon  E(\xi)\fl{\approx}{}E(\xi')$ is a 
$(SO(n),G)$-equivariant diffeomorphism
over the identity of $S^n$, then, by choosing $a':=h(a)$ and
$b':=h(b)$, one has $(\alpha',\beta')=(\alpha,\beta)$. The proof
of Lemma \ref{homos} is then complete. 
\cqfd
\end{ccote}

\begin{ccote}\label{isolift} {\it Isotropic liftings.\/}
Let $I:=[-1,1]$ and $c\colon I\to S^n$ be the parametrisation 
of the meridian arc
$c(t)=(0,\dots,\cos(\pi t/2),\sin(\pi t/2))$. 
Let $\tilde c\colon  I\to E=E(\xi)$ be a (smooth) lifting of $c$. 
As $c(t)$ is fixed
by $SO(n-1)$, one has $B\cdot \tilde c(t) = \tilde c(t)\cdot \alpha_t(B)$,
for 
$B\in SO(n-1)$. As in the proof of Lemma \ref{homos}, one checks that
this gives a smooth path $\alpha_t$ ($t\in I$) of 
homomorphisms from $SO(n-1)$ to $G$, 
which depends on the lifting $\tilde c$. 
Call $\tilde c$ {\it isotropic} 
if $\alpha_t$ is constant: 
$\alpha_t(B)=\alpha(B)$ for all $B\in SO(n-1)$.

\begin{Lemma}\label{tildec}
Any $(n,G)$-bundle admits an isotropic lifting. 
\end{Lemma}

We shall make use of connections on $(n,G)$-bundles which are  
$SO(n)$-invariant. These can be obtained by averaging any
connection (see \cite[p. 522]{BH}), since the space 
of connections is an affine space.
If a curve $u(t)$ in $E(\xi)$ is horizontal for a 
$SO(n)$-invariant connection, then $u(t)\cdot g$ and 
$A\cdot u(t)$ are horizontal. Lemma \ref{tildec} then follows
from the following

\begin{Lemma}\label{tildeconn}
Let $\xi$ be a $(n,G)$-bundle endowed with a $SO(n)$-invariant connection.
Then, any lifting $\tilde c$ of $c$ which is horizontal
is isotropic.
\end{Lemma}

\preu If $\tilde c$ is an horizontal lifting , so are 
$B\cdot\tilde c$ and $\tilde c\cdot \alpha(B)$ 
for $B\in SO(n-1)$. As 
$B\cdot\tilde c(-1)=\tilde c(-1)\cdot \alpha(B)$, one 
has $B\cdot\tilde c(t)=\tilde c(t)\cdot \alpha(B)$
for all $t\in I$.
\cqfd
\end{ccote}

\begin{ccote}\label{xialph} {\it The $(n,G)$-bundles $\xi_{\alpha,\beta}$.\/}
If $X$ is a topological space, the unreduced suspension 
$\Sigma X$ is
$$\Sigma X:= I\times X \big/ \{(-1,x)\sim (-1,x') \hbox{ and }
(1,x)\sim (1,x'), \ \forall\ x,x'\in X\}.$$
We donote by $C_\smallminus X$ the image of $[-1,1)\times X$ in $\Sigma X$ 
and by $C_\smallplus X$ those of $(-1,1]\times X$.

Let $(\alpha,\beta)$ be a pair of smooth
homomorphisms from $SO(n)$ to $G$.
Define the space $\widehat E_{\alpha,\beta}$ by
$$\begin{array}{ccl}
\widehat E_{\alpha,\beta}&:=& I\times 
SO(n)\times G\,\big /\, 
\{ (-1,A',g)\sim (-1,A,\alpha(A^\mun A')\,  g) 
\hbox{ and }\\ &&\kern 30 mm
(1,A',g)\sim (1,A,\beta(A^\mun A')\,  g),\, 
\forall A\in SO(n)\}.\end{array}$$
The space $\widehat E_{\alpha,\beta}$ admits an obvious 
free right action of $G$ and
a map $p\colon \widehat E_{\alpha,\beta}\to\Sigma SO(n)$. 
This makes a principal $G$-bundle
over $\Sigma SO(n)$; indeed, trivializations 
on $C_\pm SO(n)$ are given by
\eqncount 
\begin{equation}
\begin{array}{lllllllllll}
\hat\varphi_\smallminus  \colon  &  [t,A,g] &\mapsto & 
([t,A],\alpha(A)\, g) & 
\hbox{ if } &-1\leq t < 1\\[2pt]
\hat\varphi_\smallplus  \colon  & [t,A,g] &\mapsto & 
([t,A],\beta(A)\, g) & \hbox{ if }& 
 -1< t \leq 1.
\end{array}\end{equation}
The change of trivializations is 
\eqncount
\begin{equation}\label{transi}
\hat\varphi_\smallminus \pcirc\hat\varphi_\smallplus ^\mun ([t,A],g) = 
([t,A],\alpha(A)\beta(A)^\mun g).
\end{equation}

Now, suppose that $\alpha\, {\scriptstyle |}_{SO(n-1)}=
\beta\, {\scriptstyle |}_{SO(n-1)}$. Form the 
space $E_{\alpha,\beta}$ as the quotient
$$E_{\alpha,\beta}:=\widehat E_{\alpha,\beta}\big/
\{[t,AB,g]\sim [t,A,\alpha(B)g],\ \forall\, B\in SO(n-1)\}.$$

Let $\varepsilon\colon  SO(n)\to S^{n-1}$ be the map which associates to
a matrix its last column; it is also the projection 
$\varepsilon \colon  SO(n)\to SO(n)/SO(n-1)\iso S^{n-1}$. 
There is an map $p\colon E_{\alpha,\beta}\to\Sigma S^{n-1}$
given by $p([t,A,g])=[t,\varepsilon(A)]$
and a free $G$-action given by $[t,A,g]\cdot g_1:=[t,A,gg_1]$
As above, we check that this defines a $G$-principal bundle
over $\Sigma S^{n-1}$; the trivializations $\hat\varphi_\pm$
descend to trivializations $\check\varphi_\pm$ over $C_\pm S^{n-1}$.

The map $[t,A]\mapsto A\cdot c(t)$ descends to a homeomorphism
$f\colon \Sigma S^{n-1}\fl{\approx}{}S^n$. 
By replacing $p$ by $f\pcirc p$,
we obtain a (topological) principal $G$-bundle 
$$\xi_{\alpha,\beta} \colon  E_{\alpha,\beta}\fl{p}{} S^n.$$
Let $S^n_\pm$ be the punctured spheres 
$S^n_\pm:=f(C_\pm S^{n-1})$. The trivializations
given by 
the compositions
\eqncount
\begin{equation}\label{eq:defvarphi}
\varphi_\pm \colon  p^\mun (S^n_\pm) \hfl{\check\varphi_\pm}{} 
C_\pm S^{n-1}\times G\ \hfl{f\times {\,\rm id}}{}\ 
S^n_\pm\times G
\end{equation}
are homeomorphisms from $p^\mun (S^n_\pm)$ onto manifolds. 
The change of trivialization is a diffeomorphism, being obtained
by conjugating that of \equationref{transi} by $f$.
Therefore, $\varphi_\pm$ produce a smooth manifold structure
on $E_{\alpha,\beta}$. The map $p$ and the $G$ action are smooth.
One checks that the map
$$\begin{array}{cccccc}
SO(n)\times \widehat E_{\alpha,\beta}\to 
\widehat E_{\alpha,\beta} & \ {\rm given\ by\ } \ &
C\cdot [t,A,g]:= [t,CA,g]
\end{array}$$
descends to a smooth $SO(n)$-action on $E_{\alpha,\beta}$
which makes $\xi_{\alpha,\beta}$ a $(n,G)$-bundle.
\end{ccote}

\begin{ccote}\label{proof1B}{\it Proof of Part (i) of Theorem B.\/}
Let $\xi$ be a $(n,G)$-bundle. By Lemma \ref{tildec} there exists
an isotropic lifting $\tilde c\colon I\to E(\xi)$ of $c$. Choosing
$a:=\tilde c(-1)$ and $b:=\tilde c(1)$ produces a representative
$(\alpha,\beta)$ of $J(\xi)$ with 
$\alpha\, {\scriptstyle |}_{SO(n-1)}=
\beta\, {\scriptstyle |}_{SO(n-1)}$.
Therefore, the image of $J$ is contained in
$\calr(n,G)\times_{(n-1)}\calr(n,G)$.

Conversely, a class $P\in\calr(n,G)\times_{(n-1)}\calr(n,G)$
is represented by a pair $(\alpha,\beta)$ with 
$\alpha\, {\scriptstyle |}_{SO(n-1)}=
\beta\, {\scriptstyle |}_{SO(n-1)}$.
Let ${\mathbf 1}$ be the identity matrix in $SO(n)$ and $e$ be the
unit element of $G$.
Computing $J(\xi_{\alpha,\beta})$ with the points 
$a:=[-1,{\mathbf 1},e]$ and $b:=[1,{\mathbf 1},e]$ in 
$E_{\alpha,\beta}$ shows that $J(\xi_{\alpha,\beta})=P$.
\qed
\end{ccote}

\section{The map $\tilde J_\gamma$}\label{S:maptH}

Let $\gamma\colon SO(n-1)\to G$ be a smooth homomorphism. 
Define a set $\calr_\gamma(n,G)$ as follows: an element of  
$\calr_\gamma(n,G)$ is represented by a pair $(\alpha,\beta)$
of smooth homomorphisms from $SO(n)$ to $G$ such that
$\alpha\, {\scriptstyle |}_{SO(n-1)}=
\beta\, {\scriptstyle |}_{SO(n-1)}=\gamma$.
Two pairs $(\alpha_1,\beta_1)$ and $(\alpha_2,\beta_2)$
represent the same element of $\calr_\gamma(n,G)$ if there is
a smooth path $\{g_t\mid t\in[-1,1]\}$
in the centralizer $Z_\gamma$ of $\gamma(SO(n-1))$ 
such that $\alpha_2(A)=g_{-1}\alpha_1(A)g_{-1}^\mun$ and
$\beta_2(A)=g_{1}\beta_1(A)g_{1}^\mun$.
There is an obvious map 
$j\colon \calr_\gamma(n,G)\to \calr(n,G)\times_{(n-1)}\calr(n,G)$.

Part (i) of Theorem B, already proven in (\ref{proof1B}),
 permits us to define a map 
$\bar J\colon \cale (n,G)\to\calr(n-1,G)$ by 
$\bar J(\xi):=\res[\alpha]=\res[\beta]$. We shall 
now compute the preimage $\bar J^\mun ([\gamma])$.

\begin{Proposition}\label{TCnew}
Let $\gamma\colon SO(n-1)\to G$ be a smooth homomorphism. Then
there exists a bijection 
$\tilde J_\gamma\colon \bar J^\mun ([\gamma])
\fl{\approx}{}\calr_\gamma(n,G)$
such that $j\pcirc\tilde J_\gamma = J$.
\end{Proposition}

The proof divides into several steps.

\begin{ccote}\label{deftiga} {\it Definition of $\tilde J_\gamma$.\/}
Let $\xi$ be a $(n,G)$-bundle with $\bar J(\xi)=[\gamma]$.
Choose, using Lemma \ref{tildec}, an isotropic lifting 
$\tilde c_0\colon I\to E(\xi)$ of $c$. As $\bar J(\xi)=[\gamma]$,
the constant path $\alpha_t^0\colon SO(n-1)\to G$
associated to $\tilde c_0$ is conjugated to $\gamma$: 
there exists $g\in G$ such that $\alpha_t^0(B)=g\gamma(B)g^\mun$.
Let $\tilde c:=\tilde c_0\cdot g$. 
Like in Equation \equationref{eqconju},
one checks that $\tilde c$ is {\it $\gamma$-isotropic},
i.e. $\alpha_t=\gamma$. 
Choosing $a:=\tilde c(-1)$ and $b:=\tilde c(1)$ then produces
a pair $(\alpha,\beta)$ of smooth homomorphisms from 
$SO(n)$ to $G$ which represents a class $\tilde J_\gamma(\xi)$
in $\calr_\gamma(n,G)$.

To see that $\tilde J_\gamma$ is thus well defined, let 
$\tilde c'$ be another $\gamma$-isotropic lifting of $c$. 
The smooth path $t\mapsto g_t\in G$ defined by
$\tilde c'(t)=\tilde c(t)\cdot g_t$ then satisfies,
for all $B\in SO(n-1)$:
$$\gamma(B)=\alpha'_t(B)=g_t^\mun\alpha_t(B)g_t
=g_t^\mun \gamma(B)g_t.$$
Therefore, $g_t\in Z_\gamma$.
One has $\alpha'(A)=g_{-1}^\mun\alpha_t(A)g_{-1}$ and
$\beta'(A)=g_{1}^\mun\beta_t(A)g_1$, for all $A\in SO(n)$, 
which proves that $\tilde J_\gamma(\xi)$ does not depend on 
the choice of a $\gamma$-isotropic lifting.

Now, if $h\colon  E(\xi)\fl{\approx}{}E(\xi')$ is a 
$(SO(n),G)$-equivariant diffeomorphism
over the identity of $S^n$ and $\tilde c\colon I\to E(\xi)$ is a
$\gamma$-isotropic lifting for $\xi$, then $c':=h\pcirc\tilde c$
is a $\gamma$-isotropic lifting for $\xi'$ giving 
$(\alpha',\beta')=(\alpha,\beta)$. This proves that
$\tilde J_\gamma$ is well defined. 
\end{ccote}

\begin{ccote}\label{pfCsur} {\it Surjectivity of $\tilde J_\gamma$.\/}
Let $(\alpha,\beta)$ represent a class $P$ in $\calr_\gamma(n,G)$.
One checks that $\tilde J_\gamma(\xi_{\alpha,\beta})=P$,
using that the path $t\mapsto [t,{\mathbf 1},e]$ 
is a $\gamma$-isotropic lifting for $\xi_{\alpha,\beta}$.
\end{ccote}

\begin{ccote}\label{pfCinj}{\it Injectivity of $\tilde J_\gamma$.\/}
Let $\gamma\colon SO(n-1)\to G$ be a smooth homomorphism. 
and let $\xi$ be a $(n,G)$-bundle with $\bar J(\xi)=[\gamma]$.
There exists $a\in E(\xi)$ with $p(a)=(0,\dots ,-1)$ and
$B\cdot a=a\cdot\gamma(B)$ for all $B\in SO(n-1)$.
Choose a $SO(n)$-invariant connection on $\xi$ and let 
$\tilde c$ an horizontal lifting of $c$ with $\tilde c(-1)=a$. 
By Lemma \ref{tildeconn}, $\tilde c$ is $\gamma$-isotropic.
If $\beta\colon SO(n)\to G$ is defined by 
$A\cdot\tilde c(1)=\tilde c(1)\cdot\beta(A)$,
then $(\alpha,\beta)$ represents $\tilde J_\gamma(\xi)$.

Consider the map $\widehat \lambda\colon \widehat E_{\alpha,\beta}\to E(\xi)$
given by 
$$\widehat \lambda([t,A,g]):= A\cdot \tilde c(t)\cdot g.$$ 
The map $\widehat \lambda$ descends
to a continuous map $\lambda\colon  E_{\alpha,\beta}\to E(\xi)$ which is 
both $SO(n)$ and $G$-equivariant and which covers the
identity of $S^n$. Therefore $\xi$ and $\xi_{\alpha,\beta}$
are isomorphic as topological $(n,G)$-bundles. What
remains to prove is that $\lambda$ is a diffeomorphism,
which is only non-trivial around the fibers $E_\pm$ above 
the north and south poles.

The connection on $\xi$ provides 
a smooth trivialization of $\xi$ restricted to the punctured 
sphere $S^n_\smallminus $ (see \ref{xialph}) in the following way. 
Consider the map
$s_\smallminus \colon p^\mun(S^n_\smallminus )\to E_\smallminus $ assigning to $z$ the end
point in $E_\smallminus $ of the horizontal path through $z$
above the meridian arc through $p(z)$. Define the 
$G$-equivariant map $\sigma_\smallminus \colon p^\mun(S^n_\smallminus )\to G$ by
$s_\smallminus (z)=a\cdot \sigma_\smallminus (z)$. The required trivialization
$\tau_\smallminus \colon p^\mun(S^n_\smallminus )\to S^n_\smallminus \times G$ is
$\tau_\smallminus (z):=(p(z),\sigma_\smallminus (z))$. 

Take the trivialization $\varphi_\smallminus $ for $\xi_{\alpha,\beta}$
defined in Equations \equationref{eq:defvarphi} of (\ref{xialph}).
As $\tilde c$ is horizontal, one has
$$\tau_\smallminus \pcirc \lambda\pcirc\varphi_\smallminus ^\mun (x,g) = (x,g).$$
This, and the same for $E_\smallplus $, prove that $\lambda$ is a diffeomorphism.
We have thus established that if $\tilde J_\gamma(\xi)$ is represented
by $(\alpha,\beta)$ then the $(n,G)$-bundle $\xi$ is isomorphic to
$\xi_{\alpha,\beta}$, which proves the injectivity of 
$\tilde J_\gamma$.

The proof of Proposition \ref{TCnew} is now complete.
\qed
\end{ccote}

\section{Proof of the main results}\label{S:pmr}

\noindent
This section contains the proofs of the results stated in \S\,\ref{intro}.

\begin{ccote}{\it Proof of Theorem B.\/}\label{pfB22}
 Part (i) has already been proven in (\ref{proof1B}). 
We shall now prove Part (ii).
Let $\alpha,\beta\colon SO(n-1)\to G$ be two smooth homomorphisms
such that $\alpha\, {\scriptstyle |}_{SO(n-1)}=
\beta\, {\scriptstyle |}_{SO(n-1)}= \gamma$. 
The pair $(\alpha,\beta)$ defines a class
$[\alpha,\beta]\in\calr_\gamma(n,G)$. 
The group $Z_\gamma\times Z_\gamma$ 
acts on $\calr_\gamma(n,G)$ by 
$$(g,h)\cdot [\alpha,\beta]:=
[g\alpha g^\mun,h\beta h^\mun].$$ 
The set $J^\mun([\alpha],[\beta])\subset \bar J^\mun([\gamma])$ 
is, via the bijection 
$\tilde J_\gamma\colon \bar J^\mun ([\gamma])
\fl{\approx}{}\calr_\gamma(n,G)$ of Proposition \ref{TCnew},
in bijection with 
an orbit of the above action.
Let ``$\sim$'' be the equivalence relation on $Z_\gamma\times Z_\gamma$
defined
by $(g_1,h_1)\sim (g_2,h_2)$ iff
$(g_1,h_1)[\alpha,\beta]=(g_2,h_2)[\alpha,\beta]$.
Let 
$$\phi\colon Z_\gamma\times Z_\gamma\to
\pi_0(Z_{\alpha})\backslash\pi_0(Z_{\gamma})/\pi_0(Z_{\beta})$$
be the map defined by $\phi(g,h):=[g^\mun h]$. Part (ii) of Theorem B then
follows from the
following lemma.
\end{ccote}

\begin{Lemma} 
$(g_1,h_1)\sim (g_2,h_2)$ iff $\phi(g_1,h_1)=\phi(g_2,h_2)$.
\end{Lemma}

\preu Suppose that $(g_1,h_1)\sim (g_2,h_2)$. This means that there exist
$s_\smallminus ,s_\smallplus \in Z_\gamma$, with $[s_\smallminus ]=[s_\smallplus ]$ in $\pi_0(Z_\gamma)$, such that
the following equality
$$(g_1\alpha g_1^\mun,h_1\beta h_1^\mun)=
(s_\smallminus g_2\alpha g_2^\mun s_\smallminus ^\mun,s_\smallplus h_2\beta h_2^\mun s_\smallplus ^\mun)$$
holds in $Z_\gamma\times Z_\gamma$. This implies that 
$$g_1=s_\smallminus \,g_2\,A  \ \hbox{ and }\ h_1=s_\smallplus \,h_2\,B$$
with $A\in Z_\alpha$ and $B\in Z_\beta$ (the centralizers of the images of
$\alpha$ and $\beta$).
Therefore $g_1^\mun h_1=A^\mun g_2^\mun s_\smallminus ^\mun s_\smallplus \,h_2\,B$, which implies 
$\phi(g_1,h_1)=\phi(g_2,h_2)$.

To prove the converse, observe that
\begin{itemize}
\item $(g,h)\sim(C\,g,C\,h)$ for $C\in Z_\gamma$.
\item $(g,h)\sim (gA,hB)$ for $A\in Z_\alpha$ and $B\in Z_\beta$.
\item $(g,h)\sim(g,uh)$ for $u$ in the identity component of $Z_\gamma$.
\end{itemize}

Suppose that $\phi(g_1,h_1)=\phi(g_2,h_2)$. This means that there are 
$A\in Z_\alpha$, $B\in Z_\beta$ and $u$ in the identity component of
$Z_\gamma$
auch that $g_1^\mun h_1 = A^\mun g_2^\mun u h_2B$ ($u$ can be put in the
middle
since the identity component of $Z_\gamma$ is a normal subgroup of
$Z_\gamma$).
One then has
$$(g_1,h_1)\sim (e,g_1^\mun h_1)=(e,A^\mun g_2^\mun u h_2B)
\sim (g_2A,u\,h_2B)\sim (g_2,h_2)\ .$$
\end{proof}

\begin{proof}[Proof of Proposition C]
Let $\xi$, $\alpha$ and $\beta$ as in the statement of Proposition C.
Let $\gamma:=\alpha_{|SO(n-1)}=\beta_{|SO(n-1)}$. Then, 
$\xi\in \bar J^\mun([\gamma])$ and, by \ref{deftiga},
one has $\tilde J(\xi)=[\alpha,\beta]$ in $\calr_\gamma(n,G)$. 
By \ref{pfCinj}, $\xi=\xi_{\alpha,\beta}$.
Therefore, 
$C(\xi^\flat) = C(\xi_{\alpha,\beta}^\flat)=
[\bar\psi(\alpha,\beta)]$,
the last equality coming from Equation \equationref{transi} of \ref{xialph}
and that $C(\xi^\flat)$ can be represented by the transition function
\cite[Th. 18.4]{St}.
\end{proof}

\begin{proof}[Proof of Theorem A]
Since $SO(1)$ is trivial, $Z_\gamma=G$ which is suppose to be connected.
Therefore, Part (i) is a particular case of Part (i) of Theorem B.
Let $c\colon I\to S^n$ parametrizing the meridian arc, as in \equationref{isolift}.
Let $\alpha,\beta\colon SO(2)\to G$ be two homomorphisms representing
$J(\xi)$. One can find $a$ and $b$ so that $\alpha$ and $\beta$ are
the isotropy representations associated to $a$ and $b$. 
As $G$ is connected, the submanifold $P_0:=p^\mun (c(I))$ of $E(\xi)$
is connected and there is a smooth
lifting $\tilde c$ of $c$ such that
$\tilde c(-1)=a$ and $\tilde c(1)=b$.
As $SO(1)$ is trivial, $\tilde c$ is isotropic. Part (ii) of 
Theorem A then follows from Proposition C.
\end{proof}

\begin{proof}[Proof of Proposition D]
Recall that any element of $\pi_1(G,e)$ can be represented by a homomorphism
(a geodesic in a maximal compact subgroup $K$ of $G$, with a $K$-bi-invariant
Riemannian metric, being a $1$-parameter subgroup \cite[Ch.\,IV, \S\,6]{He}). 
Therefore, if $\eta$ is a $G$-bundle over $S^2$, there exists a
homomorphism $\alpha\colon  S^1\to G$ such that $C(\eta)=[\alpha]$. 
For $q\in\bbn$, let $\alpha_q\colon S^1\to G$
given by $\alpha_q(z):=\alpha(z)^q$. If $\alpha$ is not trivial,
the classes $[\alpha_q]$ are all distinct in $\calr(2,G)$. Indeed, 
the set $\calr(2,G)$ is in bijection
with lattice points in a Weyl chamber of the Lie algebra 
of a maximal torus of $G$ and the point representing $\alpha_q$ is
$q$ times those representing $\alpha$.

Suppose first that $\eta$ is not trivial. Hence, $\alpha$ is not trivial
and $[\alpha_{q+1},\alpha_q]$ are all different classes in 
$\calr_\gamma(2,G)$ with $[\bar\psi(\alpha_{q+1},\alpha_q])=C(\eta)$.
The result then follows from Propositions C and \ref{TCnew}.

When $\eta$ is trivial, one takes any non-trivial homomorphism 
$\alpha \colon SO(2)\to G$. The classes $[\alpha_q,\alpha_q]$ 
in $\calr_\gamma(2,G)$ represent infinitely many distinct
$SO(2)$-equivariant $G$-bundles $\xi_q$ with trivial $\xi_q^\flat$. 
\end{proof}

\begin{proof}[Proof of Proposition E]
If $n\geq 3$, the group $SO(n)$ is semi-simple and
the set $\calr(n,G)$ is finite. The latter follows from the following known
results:
\begin{itemize}
\item a homomorphism is 
determined by its tangent map at the identity (as a homomorphism
of Lie algebras).
\item the Lie algebra of $G$ contains only finitely many
semi-simple Lie subalgebras, up to inner automorphism 
\cite[Prop. 12.1]{Ri}.
\item there are only finitely many homomorphisms 
between two semisimple Lie algebras.
\end{itemize}

Also, if $G$ is compact, the group $Z_\gamma$ is compact and then
$\pi_0(Z_\gamma)$ is finite. Proposition E then follows from
Theorem B.
\end{proof} 

\begin{Remark} To remove  the hypothesis
``$G$ compact'' from Proposition E, it is enough to consider the case $G$ connected.
Indeed, $\calr_\gamma(n,G)$ is a quotient of 
$\calr_\gamma(n,G_e)$, where $G_e$ is the connected component of $e$.
One would then need the following kind of result: 
if $H$ is a compact Lie subgroup of a connected Lie group $G$, then
$\pi_0(Z(H))$ is finite. We do not know whether this true. 
\end{Remark}

\section{$SO(n+1)$-equivariant bundles}\label{S:son+1}

In this section, we describe the $(n,G)$-bundles which are 
$SO(n+1)$-equivariant $G$-bundles, for the natural action of 
$SO(n+1)$ on $S^n$. 
Let $\delta\colon  SO(n)\hfl{\approx}{} SO(n)$ 
be the conjugation by the 
diagonal $(n\times n)$-matrix  
${\rm Diag\,}(1,\dots ,1,-1)$ (or, equivalently, ${\rm Diag\,}(-1,\dots
,-1,1)$). 
If $\alpha\colon SO(n)\to G$ is a smooth homomorphism, observe that
$\res[\alpha]=\res[\alpha\pcirc\delta]$ in $\calr(n-1,G)$.

\begin{Theorem}\label{biginvar} Let $\xi$ be a $(n,G)$-bundle. 
 If $\xi$ comes from
an $SO(n+1)$-equivariant $G$-bundle then 
$J(\xi)$ is of the form $([\alpha],[\alpha\pcirc\delta])$.

 For any $[\alpha]\in \calr(n-1,G)$ there is a unique $\xi\in \calr(n,G)$
which comes from a $SO(n+1)$-equivariant $G$-bundle and 
such that $J(\xi)=([\alpha],[\alpha\pcirc\delta])$.
\end{Theorem}

\preu For $\theta\in [0,\pi]$, let $R_\theta\in SO(n+1)$ be the rotation
of angle $\theta$ in the plane of the last $2$ coordinates. 
Let $R:=R_\pi$, the diagonal matrix with entries
$(1,\dots ,1,-1,-1)$. 

Let $\xi$ be an $SO(n+1)$-equivariant bundle.
Choose $a,b\in E(\xi)$, with
$p(a)=(0,\dots ,-1)$ and let $b:=R\cdot a$. For 
$A\in SO(n)$, one has $R^{-1}AR=\delta (A)$ and
\eqncount 
\begin{equation}\label{biginvarpfeq1}
\begin{array}{llllll}
b\,\beta(A) &=& A\cdot b = A\cdot(R\cdot a) = R\cdot(R^{-1}AR)\cdot a \\[3pt]
&=& R\cdot a\,\alpha(\delta (A)) = b\,\alpha(\delta (A)),
\end{array}\end{equation}
whence $\beta=\alpha\pcirc\delta$, which proves Part (i).

Let $\alpha\colon SO(n)\to G$ be a smooth homomorphism and set
$\gamma:=\alpha_{|SO(n-1)}$. 
Suppose that $\xi$ is an $SO(n+1)$-equivariant $G$-bundle with $a\in
E(\xi)$ such
that $p(a)=(0,\dots ,-1)$ and $A\cdot a=a\,\alpha(A)$ for $A\in SO(n)$. 
Then $\xi\in \bar J^{-1}([\gamma])$. let $c\colon I\to E(\xi)$ be the curve
$c(t):=R_{\theta(t)}\cdot a$, where $\theta(t):=\frac{\pi}{2}(t-1)$. Using
that
$R_\theta^{-1}BR_\theta = B$ for $B\in SO(n-1)$, on checks, as in equation 
\equationref{biginvarpfeq1}, that $c$ is a $\gamma$-isotropic lifting. By 
\equationref{deftiga} and equation \equationref{biginvarpfeq1}, 
one has $\tilde J_\gamma (\xi)=([\alpha],[\alpha\pcirc\delta])$ in
$\calr_\gamma(n,G)$.
By Proposition \ref{TCnew}, this proves that $\xi$ is determined by $\alpha$,
which proves the uniqueness statment of Part (ii).

It remains to construct, for a smooth homomorphism $\alpha \colon  SO(n)\to G$,
an $SO(n+1)$-equivariant $G$-bundle with
$J(\xi)=([\alpha],[\alpha\pcirc\delta])$.
Consider the
map $p\colon SO(n+1)\to S^n$ sending a matrix to its last column. This
makes an $SO(n+1)$-equivariant $SO(n)$-bundle (the principal 
$SO(n)$-bundle associated to the tangent bundle of $S^n$). Let $\xi$ 
be the $G$-bundle obtained by the Borel construction, using
the homomorphism $\alpha\pcirc\delta\colon SO(n)\to G$
\eqncount
\begin{equation}\label{conborel}
E(\xi):= SO(n+1)\times G\big/ \{(BA,g)=(B,\alpha(\delta(A))g)\}.
\end{equation}
This is an $SO(n+1)$-equivariant $G$-bundle. 
Choosing $a:=(R,e)$ and $b:=(I,e)$,
one sees that 
$$ A\cdot a = (AR,e) = (R\delta(A),e) = (R,\alpha(A)e) = a\,\alpha(A)$$
and
$$ A\cdot b = (A,e) = (I,\alpha(\delta(A))e) = b\,\alpha(\delta(A)).$$
Therefore, $\xi$ is an $SO(n+1)$-equivariant $G$-bundle with
$J(\xi)=([\alpha],[\alpha\pcirc\delta])$. 
\cqfd

\begin{Remark}  Theorem \ref{biginvar} and its proof show that
any $SO(n+1)$-equivariant $G$ bundle is derived from the tangent bundle
to $S^n$ by the Borel construction (formula \equationref{conborel}). This can be 
compared with \cite[\S\ 6]{Ha}.
\end{Remark}

\section{Examples and applications}\label{exap}

\noindent\it Notation: \rm if $X$ is a set, we denote by 
$\Delta X$ the diagonal in $X\times X$.

\begin{ccote}\label{2Um} {\it $n=2$ and $G=U(m)$.\/}
A homomorphism $\alpha \colon  S^1\to U(m)$ has, up to conjugacy,
a unique diagonal form 
$\alpha(z)={\rm Diag\,}(z^{p_1},\dots ,z^{p_m})$, with
$p_1\geq\cdots\geq p_m$. The same holds for $\beta$.
By Theorem A, $\cale(2,U(m))$ is then in bijection
with the set of pairs $(p,q)$ of $m$-tuples of 
non-decreasing integers. In $\pi_1(U(m))=\bbz$, one
has $[\alpha]=\sum_{i=1}^mp_i$ and $[\beta]=\sum_{i=1}^mq_i$, so,
by Proposition C:
$$C(\xi^\flat)=\sum_{i=1}^m(p_i-q_i).$$
If one wishes instead to characterize $\xi^\flat$ by its first
Chern number $c(\xi^\flat)\in H^2(S^2)=\bbz$, then 
$c(\xi^\flat)=-C(\xi^\flat)$ \cite[p. 445]{Mi}.

For instance, if $\tau$ is the unit tangent bundle over $S^2$ with
the natural action, then $\alpha(z)=z$, $\beta(z)=z^{-1}$,
so $C(\tau)=-2$ and $c(\tau)=2=\chi(S^2)$. 

Note that, if $\xi$ comes from a $SO(3)$-bundle, then, by Theorem
\ref{biginvar},
one has $q_i=-p_i$ (like for $\tau$ above). In particular $c(\xi^\flat)$
must be even
(see also \cite[(6.3)]{Ha}).
\end{ccote}

\begin{ccote}\label{2O2} {\it $n=2$ and $G=O(2)$.\/}
For $q\in\bbz$, let $\alpha_q\colon SO(2)\to O(2)$ be the homomorphism $A\to
A^q$. The set
$\calr(2,O(2))$ is in bijection with $\bbn$ 
given by $\alpha_q\mapsto |q|$. The same recipe produces a bijection
$\pi_1(O(2),e)/\pi_0(O(2))\iso\bbn$. 
Let $\xi_{p,q}:=\xi_{\alpha_p,\alpha_q}$ be the $(2,O(2))$-bundles
constructed in \equationref{xialph}. By Proposition \ref{TCnew} and \equationref{pfCsur}, 
each $\cale(2,O(2))$ is represented by some $\xi_{p,q}$, with the only
relation $\xi_{p,q}=\xi_{-p,-q}$. One has $J(\xi_{p,q})=(|p|,|q|)$
and $C(\xi_{p,q})=|p-q|$. Therefore, $J^\mun(r,s)$ contains 1 element
if $rs=0$ and 2 otherwise. 
\end{ccote}

\begin{ccote}\label{2SOm} {\it $n=2$ and $G=SO(m)$ with $m\geq 3$.\/}
A maximal torus of $SO(m)$ is formed by matrices containing 2-blocks
concentrated
around the diagonal, so isomorphic to $SO(2)^k$ where $k=[m/2]$. As in
\ref{2Um},
by Theorem A, $\cale(2,SO(m))$ is then in bijection
with the set of pairs $(p,q)$ of $k$-tuples of 
non-decreasing integers. The bundle $\xi^\flat$ is determined
by its second Stiefel-Whitney number $w(\xi^\flat)\in\bbz_2$ which is then
given by
$$w(\xi^\flat)=\sum_{i=1}^{[\frac{m}{2}]}(p_i-q_i) \qquad ({\rm mod\,} 2).$$
Again, $\xi$ comes from a $SO(3)$-equivariant bundle iff $q_i=-p_i$ and then
$\xi^\flat$ is trivial.
\end{ccote}

\begin{ccote}\label{2k+1} {\it $n=2k+1\geq 3$ and $G$ is a compact classical
group other than $SO(2m)$.\/}
The important thing is that $SO(n-1)$ contains a
maximal torus of $SO(n)$. Therefore, by \cite[Ch. 6, Corollary 2.8]{BT},
for any embedding $\psi\colon G\hookrightarrow U(m)$, the representations
 $\psi\pcirc\alpha$ and $\psi\pcirc\beta$ 
are conjugate in  in
$U(M)$. For  
$G$ a compact classical group other than $SO(2m)$, this implies that 
$\alpha$ and $\beta$ are conjugate in $G$ \cite[Pro. 8, p.56]{On}. Therefore, 
the image of $J$ is the diagonal $\Delta\calr(n,G)$. 
We do not know whether this true for $G=SO(2m)$ (see \cite[Rm. p.57]{On} for
a possible source of counter-examples).
\end{ccote}

\begin{ccote}\label{2kSOn} {\it $n=2k+1\geq 3$ and $G=SO(n)$.\/}
The set $\calr(n,SO(n))$ has just two elements,
represented by the trivial homomorphism and the
identity $\rm id$ of $SO(n)$.
If $\iota\colon  SO(n-1)\subset SO(n)$ denotes the inclusion, then
$Z_{[\iota]}$ contains 2 elements, represented by the identity matrix 
and the diagonal matrix $D$ with entries 
$(-1,\dots ,-1,1)$. Let $\delta
$ be the inner automorphism of 
$SO(n)$ given by the conjugation by $D$. By Theorem B,
the set $\cale (n,SO(n))$ for $n=2k+1\geq 3$
then contains 3 elements:
\begin{enumerate}
\renewcommand{\labelenumi}{(\roman{enumi})}
\item  the trivial bundle $S^n\times SO(n)$ with the action
$A\cdot (z,B)=(A\cdot z,B)$. 
The isotropy representations are both trivial.
\item  the trivial bundle $S^n\times SO(n)$ 
with the action
$A\cdot (z,B)=(A\cdot z,AB)$. It is characterized by 
$\bar J(\xi)=[\iota]$ and 
$J_\iota(\xi)=([{\rm id}],[{\rm id}])$.
This does not come from an $SO(n+1)$-equivariant bundle.
\item the principal $SO(n)$-bundle $\calt S^{n}$ associated with
the tangent bundle of $S^n$. It is characterized by 
$\bar J (\calt S^{n})=[\iota]$ and 
$J_\iota(\calt S^{n})=([{\rm id}],[\delta])$.
This comes from an $SO(n+1)$-equivariant bundle.
\end{enumerate}

Observe that $([{\rm id}],[\delta])=([\delta],[{\rm id}])$ in 
$\calr_\iota (n,SO(n))$. By Proposition C, this implies
that $C(\calt S^{n})=-C(\calt S^{n}
)$. This proves again the classical fact
that $C(\calt S^{2k+1})\in\pi_{2k}(SO(2k+1))$ 
is of order 2 \cite[Cor. IV.1.11]{Br}. 
\end{ccote}

\begin{ccote}\label{2kSOn6} {\it $n=2k\geq 6$ and $G=SO(n)$.\/}
The set $\calr(n,SO(n))$ contains 3 elements,
represented by the trivial homomorphism, the
identity $\rm id$ of $SO(n)$
and the conjugation $\delta$ by the diagonal matrix with entries 
$(-1,\dots ,-1,1)$.
The non-trivial homomorphisms restrict to the inclusion
$\iota\colon  SO(n-1)\subset SO(n)$. The group
$Z_{\iota}$ is trivial. Therefore,
$J$ is injective and the set $\cale (n,SO(n))$ for $n=2k$
then contains 5 elements:
\begin{enumerate}
\renewcommand{\labelenumi}{(\roman{enumi})}
\item  the trivial bundles $S^n\times SO(n)$ with the actions
$A\cdot (z,B)=(A\cdot z,B)$,  
$A\cdot (z,B)=(A\cdot z,AB)$ and  
$A\cdot (z,B)=(A\cdot z,\delta(A)B)$. Their images
by $J$ is the diagonal $\Delta\calr(n,SO(n))$.
\item the principal $SO(n)$-bundle $\calt S^{n}$ associated with
the tangent bundle of $S^n$. 
One has $J(\calt S^{n})=([{\rm id}],[\delta])$.
\item the $(n,G)$-bundle $-\calt S^{n}$ with
$J(-\calt S^{n})=([\delta],[{\rm id}])$. Its underlying
$SO(n)$-principal bundle is stably trival 
with Euler number $-2$. 
\end{enumerate}

The trivial bundle with action $A\cdot (z,B)=(A\cdot z,B)$, as well
as the bundles in (ii) and (iii) are the ones coming from 
$SO(n+1)$-equivariant bundles.
\end{ccote}

\begin{ccote}\label{43} {\it $n=4$ and $G=SO(3)$.\/}
The groups $SO(4)$ and $SO(3)$ are built up out of the unit quaternions 
$S^3$ by 
$SO(4)\iso (S^3\times S^3)\big/\{(1,1),(-1,-1)\}$ and
 $SO(3)\iso S^3/\{\pm 1\}$. Recall that these isomorphisms
are constructed as follows:  the orthogonal transformation $A_{p,q}\in SO(4)$
associated to $(p,q)\in S^3\times S^3$ is $A_{p,q}(x):=px\bar q$,
where $x\in\bbh$ is a quaternion and $\bbh$ is made isomorpohic to 
$\bbr^4$ by choosing $(i,j,k,1)$ as a basis.
The correspondence $p\to A_{p,p}$ then induces 
the inclusion $\iota\colon SO(3)\subset SO(4)$. As for the automorphism
$\delta$ of $SO(4)$ of \S\ \ref{S:son+1}, the conjugation by 
$D:={\rm Diag\,}(-1,\dots ,-1,1)$, as $Dx=\bar x$, one checks easily
that $\delta(A_{p,q})=A_{q,p}$.

The non-equivariant isomorphism class of a 
$SO(3)$-principal bundle $\eta$ is characterized by 
$C(\eta)\in \pi_3(SO(3))=\pi_3(S^3)=\bbz$. It is also
determined by its first Pontrjagin number 
$p(\eta)\in 4\,\bbz$, with the relation
$p(\eta)=4\,C(\eta)$.

The set $\calr (4,SO(3))$ contains 
3 elements represented by
the trivial homomorphism and those induced by the 
projections $S^3\times S^3 \to S^3$ given by
$\sigma_1(p,q):=p$, $\sigma_2(p,q):=q$.
The last two restrict over $SO(3)$ to the identity ${\rm id}$ 
of $SO(3)$. The group
$Z_{\iota}$ being trivial, the map $J$ is injective. 
This shows that the set $\cale (4,SO(3))$
contains 5 elements:
\begin{enumerate}
\renewcommand{\labelenumi}{(\roman{enumi})}
\item  the trivial bundles $S^4\times SO(3)$ with the actions
$A\cdot (z,B)=(A\cdot z,B)$ and  
$A\cdot (z,B)=(A\cdot z,\sigma_i(A)B)$ for $i=1,2$. 
Their images by $J$ is the 
diagonal $\Delta\calr(4,SO(3))$.
\item the principal $SO(3)$-bundle $\calh \colon  \bbr P^7\to S^4$ 
coming from the quaternionic Hopf bundle $S^7\to S^4$;
the $SO(4)$-action comes from the $SU(2)\times SU(2)$-action
on $S^7$ given by $(p,q)\cdot (z_1,z_2)=(pz_1,qz_2)$.
One has $J(\calh)=([{\sigma_1}],[\sigma_2])$ and $p(\calh^\flat)=4$.
\item the $(n,G)$-bundle $-\calh$ with
$J(-\calh)=([{\sigma_2}],[\sigma_1])$ and $p(\calh^\flat)=-4$. 
\end{enumerate} 

The trivial bundle with action $A\cdot (z,B)=(A\cdot z,B)$, as well
as the bundles in (ii) and (iii) are the ones coming from 
$SO(5)$-equivariant bundles.
\end{ccote}

\begin{ccote}\label{44} {\it $n=4$ and $G=SO(4)$.\/}
Taking the notations of \ref{43}, the set
$\calr (4,SO(4))$ contains 5 elements represented by:
\begin{itemize}
\item the trivial homomorphism.
\item those induced by $\sigma_1(p,q):=(p,p)$ and $\sigma_2(p,q):=(q,q)$.
\item the identity ${\rm id}$ of $SO(4)$.
\item the homomorphism $\delta(p,q):=(q,p)$.
\end{itemize}
The non-trivial homomorphisms all restrict to 
to $\iota$ over $SO(3)$. The group
$Z_\iota$ is trivial and then $J$ is injective. 
The non-equivariant isomorphism class of a 
$SO(4)$-principal bundle $\eta$ is characterized by 
$C(\eta)\in\pi_3(SO(4))=\pi_3(S^3\times S^3)=\bbz\oplus\bbz$.
More usually, one takes the pair of integers $(p(\eta),e(\eta))$
formed by the first Pontrjagin number ($\in 2\bbz$) and the Euler number of
$\eta$. 
The linear map which sends $C(\eta)$ to $(p(\eta),e(\eta))$
has matrix $\big({}^2_1\,{}^{{\hphantom{-}}2}_{-1}\big)$ (as can be checked
on the examples
below).

One sees that the set $\cale (4,SO(3))$
then contains 17 elements:
\begin{enumerate}
\renewcommand{\labelenumi}{(\roman{enumi})}
\item  the trivial bundles $S^4\times SO(3)$ with the 5 
actions using the above representations.
Their images by $J$ are just the diagonal elements $\Delta\calr(4,SO(4))$.
\item the principal $SO(4)$-bundle 
$\widehat\calh$ whose total space is
$E(\widehat\calh):=\bbr P^7\times_{SO(3)}SO(4)$. 
One has $J(\widehat\calh)=([{\sigma_1}],[\sigma_2])$ 
and $C(\widehat\calh^\flat)=(1,1)$ 
therefore its characteristic classes are $(p(\calh^\flat),e(\calh^\flat))=(4,0)$.
Also its ``inverse" $-\widehat\calh$, with 
$J(-\widehat\calh)=([{\sigma_2}],[\sigma_1])$ 
and $(p(\calh^\flat),e(\calh^\flat))=(-4,0)$.
\item the principal $SO(4)$-bundle $\calt S^{4}$ associated with
the tangent bundle of $S^4$. 
One has $J(\calt S^{4})=([{\rm id}],[\delta])$. Therefore,
$C(\calt S^{4})=(1,-1)$ and 
$(p(\calt^\flat),e(\calt^\flat))=(0,2)$. 
Again, one can consider its inverse.
\item the $(n,G)$-bundles $\xi_{i}$ ($i=1,2$) with
$J(\xi_{i})=([{\rm id}],[\sigma_i])$ and their inverses
$-\xi_{i}$. They satisfy 
$C(\xi_1)=(0,-1)$ and $C(\xi_2)=(1,0)$, or, equivalently:
$$(p(\xi_1^\flat),e(\xi_1^\flat))=(-2,1) \ \hbox{ and }\
(p(\xi_2^\flat),e(\xi_2^\flat))=(2,1).$$

\item the $(n,G)$-bundle $\xi_{i,\delta}$ ($i=1,2$) with
$J(\xi_{i})=([{\delta}],[\sigma_i])$ and their inverses.
They satisfy
$C(\xi_{1,\delta})=(-1,0)$ and
$C(\xi_{1,\delta})=(0,1)$,
or, equivalently:
$$(p(\xi_{1,\delta}^\flat),e(\xi_{1,\delta}^\flat))=(-2,-1) \ \hbox{ and }\
(p(\xi_{2,\delta}^\flat),e(\xi_{2,\delta}^\flat))=(2,-1).$$
\end{enumerate} 

Only the trivial bundle with action $A\cdot (z,B)=(A\cdot z,B)$ and
the bundles in (ii) and (iii) come from 
$SO(5)$-equivariant bundles.

\end{ccote}

\begin{ccote}\label{nUm} {\it $G=U(m)$.\/}
In order to have non-trivial $(n,U(m))$-bundles, one must have 
${\rm dim\,}U(m)> {\rm dim\,}SO(n)$. We check that we 
are then in the stable range, where, by Bott periodicity,
$$\pi_{n-1}(U(m))\approx\pi_{n-1}(U(m+k))\approx
\left\{\begin{array}{lllll}
0 & \hbox{if $n$ is odd}\\
\bbz & \hbox{if $n$ is even.}\end{array}\right.$$
Problem: {\it which integers occur as $C(\xi^\flat)$ for
a $(n,U(m))$-bundle $\xi$\, ?}
\end{ccote}

\section{A more general setting}\label{genset}

The orthogonal action of $SO(n)$ on $S^n$ is an example of
the {\it special} $\Pi$-manifolds defined by J\"anich \cite[1.2]{Ja}.
Other examples include the ``cohomogeneity one" actions studied
by E.~Straume \cite{Straume} (see \cite{GZ}
for a recent application and other
references). In this section we give the classification of equivariant
$(\Pi, G)$-bundles over special $\Pi$-manifolds. 
We will assume in this section
 that $\Pi$ and $G$ are both {\it compact\/} Lie groups.

Let $X$ be a smooth, connected, closed $n$-dimensional manifold with a
smooth $\Pi$-action. Choose a $\Pi$-invariant Riemannian metric
on $X$,  and then each tangent space $T_x X$ contains a $\Pi_x$
invariant subspace $V_x$ perpendicular to the orbit $\Pi\cdot x$.
Then $X$ is called a {\it special} $\Pi$-manifold if for each
$x \in X$ the representation of $\Pi_x$ on the normal space
$V_x$ is the direct sum of a trivial representation and a transitive
representation.

It follows that the orbit space $M = X/\Pi$ admits a natural structure
as a topological manifold with boundary
\cite[1.3]{Ja}, with dimension equal to $n - \dim (\Pi/H)$
where $H$ is the principal isotropy type. 
Under the ``functional" smooth structure \cite[VI.6]{B}, the orbit
space $M$ is a smooth manifold with boundary.
The pair $(X, \pi\colon X \to M)$
is called a special $\Pi$-manifold over $M$.

Special $\Pi$ manifolds over $M$ were classified by J\"anich
\cite[3.2]{Ja}, and independently  by W.-C. Hsiang and W.-Y. Hsiang
\cite{Hs} (see also \cite[V.5, VI.6]{B}).

 Let
$\partial M_A =\{B_\alpha\}_{\alpha \in A}$ denote the set of boundary
components of $M$.
An {\it admissible\/} isotropy group system $(H, U_A)$ over $M$ consists
of
a closed subgroup $H$ of $\Pi$ and a set $U_A = \{U_\alpha\}_{\alpha \in
A}$
of closed subgroups in $\Pi$ containing $H$, with the property that
for each $\alpha \in A$ there exists a transitive representation in which
$H$ appears as the isotropy group of a non-zero vector. Let $\Gamma =
N(H)/H$
and 
$\Omega_\alpha = N(U_\alpha) \cap N(H)/H$ for each $\alpha \in A$.
The idea of the classification is to re-construct $X$ from the 
unique principal
$\Gamma$-bundle $P$ over $M$ such that
$P_{|M_0}= \{ x \in X \vv \Pi_x = H\}$, 
and  a reduction of the structural group
of $P|B_\alpha$ to $\Omega_\alpha$ over each of
the boundary components $B_\alpha$. 


Let $B_\alpha \times [0,1]$ be
a collar neighbourhood of some component $B_\alpha$ in $M$, and 
let $Y_\alpha = \pi^{-1}(B_\alpha\times [0,1])$
 denote its pre-image in $X$. The key
fact is the following identification of $Y_\alpha$ as a $\Pi$-space.
 
\begin{Theorem}{\label{tube}} 
{\rm (Tube Theorem \cite[V.4.2]{B})} 
Let $Y=Y_\alpha$,
$B = B_\alpha$ and $\Omega = \Omega_\alpha$. There
exists  
a right $\Omega$-principal bundle $Q = Q_\alpha$ over $B$,
and a $\Pi$-equivariant diffeomorphism
$$M_\psi \times_{\Omega} Q \stackrel{\approx}{\longrightarrow} Y$$
commuting with the projection to $[0,1]$, where $M_\psi$ denotes
the mapping cylinder of the canonical projection 
$\psi\colon \Pi/H \to \Pi/U$.
\end{Theorem} 

Let $X_0 = X - \bigcup_\alpha \pi^{-1}(B_\alpha \times [0,1/2))$
and $M_0 = X_0/\Pi$.
The Tube Theorem shows that $X$ is $\Pi$-diffeomorphic to the union  
$$ X = X_0 \cup \bigcup_\alpha Y_\alpha
 = \Pi/H \times_{\Gamma} P \cup \bigcup_\alpha
 M_{\psi_\alpha} \times_{\Omega_\alpha} Q_\alpha $$
with the identification on the overlaps $B_\alpha \times [1/2,1]$
induced  by a reduction of structural groups
$P | B_\alpha \cong \Gamma \times_{\Omega_\alpha} Q_\alpha$.

\medskip
Two special $\Pi$-manifolds $(X_1, \pi_1)$ and $(X_2, \pi_2)$
over $M$ are called {\it isomorphic} when there exists a 
$\Pi$-equivariant diffeomorphism $f\colon X_1 \to X_2$ so that
the induced diffeomorphism $\bar f \colon M \to M$ is the identity.
By the the smooth isotopy lifting theorem of G. Schwarz 
\cite[Corollary 2.4]{Sch} this is equivalent to J\"anich's original
definition where $f$ was assumed to be only isotopic
to the identity, by an isotopy fixing $\partial M$ pointwise.
Let ${\mathcal S}[H, U_A]$ denote the set of isomorphism classes of special 
$\Pi$-manifolds over $M$, with isotropy group system
fine equivalent to $(H, U_A)$ (cf. \cite[\S 2]{Ja}).

\begin{Theorem}
{\label{special}} {\rm (\cite[3.2]{Ja})}
 Let $(H, U_A)$ be an
admissible isotropy group system over $M$, where $M$ is a smooth, compact
connected manifold with boundary. Then
$${\mathcal S}[H, U_A] \cong [M, \partial M_A; B\Gamma, B\Omega_A]\ .$$
\end{Theorem}

In order to classify equivariant $(\Pi, G)$-bundles $(E, p)$ over a
special
$\Pi$-manifold $X$, called {\it special $(\Pi, G)$-bundles\/} for short,
 we will generalize the results of Lashof \cite{L3} to describe
the bundles over $\Pi$-spaces $X_0$ and $Y_\alpha$ with one orbit,
and then follow J\"anich's method to glue the pieces together.

First some general definitions: if $\rho\colon H \to G$ is a (smooth)
homomorphism,  $[\rho]$ denotes the set of homomorphisms
$\rho'\colon H \to G$ such that $\rho'(h) = g\rho(h)g^{-1}$ for 
some $g \in G$ and all $h\in H$. 
We will say that {\it the fibre over 
$x$ belongs to $[\rho]$\/} if for each $z\in
p^{-1}(x)$ there exists $\rho' \in [\rho]$ such that 
$hz = z \cdot \rho'(h)$ for all $h\in H$. Then let

$$X^{[\rho]} = \{ x \in X^H \vv {\rm \ the\ fibre\ over\ }
x {\rm \  belongs\ to\ } [\rho] \}$$
Let $X^{[\rho]}_0 = X_0 \cap X^{[\rho]}$ and 
notice that
$X^{[\rho]}_0 \subset P = \{x \in X \vv \Pi_x = H\}$.
 By \cite[Lemma 1.2]{L3}, the space
$$E^\rho = \{ z \in E \vv hz = z\cdot \rho(h), \ \forall h \in H\}$$
is an $Z_\rho$-bundle over $X^{[\rho]}$,
 where $Z_\rho$ is the centralizer
of $\rho$ in $G$. The group $\widehat G = \Pi \times G$ has a left
action on the
total space $E$ given
by the formula $(\gamma, g)\cdot z = \gamma z\cdot g^{-1}$ for
any $(\gamma, g) \in \widehat G$ and any $z \in E$. 
Let us set
$$H\langle \rho \rangle := \{ (h, \rho(h)) \vv h \in H\} \subset \Pi \times
G$$
and 
$$\Gamma\langle \rho\rangle :=
 N(H\langle \rho \rangle)/H\langle \rho \rangle\ .$$
Then $E^\rho$ is just the fixed set of $H\langle\rho\rangle$ in $E$ under
this action.

Two special $(\Pi, G)$-bundles $(E_1, p_1)$ and $(E_2, p_2)$
over $M$ are called {\it equivalent}
when there exists a $\Pi$-equivariant $G$-bundle isomorphism
$\phi\colon E_1 \to E_2$, covering a $\Pi$-equivariant diffeomorphism
$f\colon X_1 \to X_2$, so that
the induced diffeomorphism $\bar f \colon M \to M$ is the identity.

We now give the classification of the
part of $(E,p)$ lying over $M_0$. After the remarks above, we see
that it follows directly from the slice theorem \cite[II.5.8]{B}.

\begin{Theorem}{\label{principal}}
 {\rm (\cite[1.9]{L3})} Let $\rho\colon H \to G$ be
a smooth homomorphism. The equivalence classes of 
 $(\Pi, G)$-equivariant
bundles $(E_0, p_0)$ over $M_0$ with all fibres belonging to $[\rho]$
are in bijection with the homotopy classes of maps 
$[M_0, B\Gamma\langle \rho\rangle]$.
\end{Theorem}

We next define an {\it admissible\/} $(\Pi, G)$-isotropy group system over
$M$ to be a set 
$(H, U_A; \rho,  \rho_A)$,
 where $H$ and $U_A$ are as
above,
$\rho\colon H \to G$ is a homomorphism, 
and $\rho_A = \{\rho_\alpha\}_{\alpha \in A}$ is a set of homomorphisms
$\rho_\alpha\colon U_\alpha \to G$ such that $\rho_\alpha | H = \rho$.
We define
$$\Omega_\alpha\langle \rho_\alpha \rangle = 
N(U_\alpha\langle \rho_\alpha \rangle) \cap
N(H\langle \rho \rangle)/H\langle \rho \rangle$$
and
 $$\Omega_A\langle \rho_A \rangle =
 \{\Omega_\alpha\langle \rho_\alpha \rangle\}_{\alpha\in A}\ .$$
A special  $(\Pi, G)$-bundles $(E, p)$ {\it realizes\/} an 
admissible $(\Pi, G)$-isotropy group system $(H, U_A; \rho, \rho_A)$
over $M$ if 
\begin{enumerate}
\renewcommand{\labelenumi}{(\roman{enumi})}
\item there exists points $\{y_\alpha \in Y_\alpha\}$ such
that $\Pi_{y_\alpha} = U_\alpha$,
\item for each $y_\alpha$, the normal space $V_{y_\alpha}$ to the orbit
$\Pi\cdot y_\alpha$ has a point $z_\alpha \in V_{y_\alpha}$ with
$\Pi_{z_\alpha} =H$,
\item the images  $c_\alpha(t)$ of rays in $V_{y_\alpha}$ joining
$z_\alpha$ to $y_\alpha$, have isotropic
liftings $\tilde c_\alpha (t)$ to $E$ such that 
$\tilde c_\alpha (t) \subset E^\rho$ for $0\leq t < 1$ and
$\tilde c_\alpha (1) \subset E^{\rho_\alpha}$.
\end{enumerate}
\noindent
It is not difficult to check (following \cite[\S 2]{Ja}) that
every special  $(\Pi, G)$-bundle $(E, p)$ over $M$ realizes some
admissible $(\Pi, G)$-isotropy group system $(H, U_A; \rho, \rho_A)$.
This isotropy group system is unique up to a 
natural notion of equivalence, 
extending the ``fine-orbit structure" of J\"anich.

We  say that two isotropy group systems
 $(H, U_A;\rho,\rho_A)$
and $(H', U'_A;\rho',\rho'_A)$
 are {\it fine-equivalent\/} if the
following
conditions hold:

\begin{enumerate}
\renewcommand{\labelenumi}{(\roman{enumi})}
\item there exists an element $\gamma \in  \Pi \times G$
such that 
$H'\langle \rho'\rangle =\gamma H\langle \rho \rangle\gamma^{-1}$, and
\item there exist $n_\alpha \in NH\langle \rho \rangle$
such that 
$U'_\alpha\langle\rho'_\alpha\rangle =
 (\gamma n_\alpha) U_\alpha\langle\rho_\alpha\rangle(\gamma n_\alpha)^{-1}$.
\end{enumerate}
Let ${\mathcal S}[H, U_A;\rho,\rho_A]$ denote the set of 
equivalence classes of special 
$(\Pi,G)$-bundles over $M$ realizing the given
$(\Pi, G)$-isotropy group system, up to fine equivalence.

\begin{Theorem}{\label{bundles}} 
 Let $(H, U_A; \rho, \rho_A)$ be an
admissible $(\Pi, G)$-isotropy group system over $M$,
 where $M$ is a smooth, compact
connected manifold with boundary. Then
$${\mathcal S}[H, U_A;\rho,\rho_A]
 \cong
 \big [M, \partial M_A;B\Gamma\langle \rho\rangle, B\Omega_A\langle \rho_A \rangle\big ]\ .$$
\end{Theorem}

\preu
Suppose that we are given an admissible
$(\Pi, G)$-isotropy group system over $M$.
Let $(E,p)$ be a special
$(\Pi,G)$-bundle over $M$ realizing the given
$(\Pi, G)$-isotropy group system.
 By restricting the bundle
to $M_0$, we get a  map $\omega_0\colon M_0 \to B\Gamma\langle
\rho\rangle$ classifying
the principal $\Gamma\langle \rho\rangle$-bundle 
$P\langle \rho\rangle$
(which completely determines
$(E_0, p_0)$) by Theorem \ref{principal}.
We can apply the Tube Theorem \ref{tube}
 to the $\widehat G =\Pi\times G$
action on $p^{-1}(Y_\alpha) \subset E$, since
$z \in E^\rho$ means that 
$\widehat G_z = H\langle \rho \rangle$ and similarly for
$z \in E^{\rho_\alpha}$.
 This identifies the restriction
of our bundle to the part over $B_\alpha \times [0,1]$ as
$M_{\varphi_\alpha} \times_{\Omega_\alpha\langle \rho_\alpha \rangle}
Q\langle
\rho_\alpha \rangle$,
where 
$$\varphi_\alpha\colon \widehat G/H\langle \rho \rangle \to 
\widehat G/U_\alpha\langle \rho_\alpha \rangle$$
is the $\widehat G$-equivariant projection, and
 $Q\langle \rho_\alpha \rangle$
is a principal right $\Omega_\alpha\langle \rho_\alpha \rangle$-bundle
over
$B_\alpha$.
The classifying map $\omega_0$ for $P\langle \rho \rangle$ 
therefore  extends to
a map 
$$\omega\colon (M, \partial M_A) \to
 (B\Gamma\langle \rho\rangle,B\Omega_A\langle \rho_A \rangle)$$
where the notation means that each boundary component $B_\alpha$
is mapped into the $\alpha$-component of 
$B\Omega_A\langle \rho_A \rangle$.
The restriction of $\omega$ to $B_\alpha$ classifies
$Q\langle \rho_\alpha \rangle$.
This shows that $(E,p)$ is determined up to equivalence by $\omega$.

Conversely, if we are given a map $\omega$ as above we can reconstruct
a special $(\Pi, G)$-equivariant bundle over $M$ realizing the isotropy
group system, up to fine equivalence.
It can be checked that this bundle is unique up to
equivalence.
\cqfd

\begin{Remark}\label{rkPiG} {Note that a special $\Pi\times G$-manifold
over $M$
is a special $(\Pi,G)$-bundle  over $M$ precisely when
 the subgroup
$1\times G$ acts freely on the total space. 
The isotropy group system
for the bundle is just the collection of isotropy groups for 
the $\Pi\times G$-action. This observation shows that Theorem
\ref{bundles}  follows from J\"anich's results.}
\end{Remark}

\begin{Corollary} Let $\Pi$ and $G$ be compact Lie groups.
 The set of special
$(\Pi, G)$-equivariant bundles over $M$
is finite, provided that $\dim M \leq 1$
 and the isotropy groups are semi-simple.
\end{Corollary}
\preu
This is proved using  \cite{Ri}, as for the finiteness of
 ${\mathcal E}(n,G)$.
\cqfd

\medskip
To conclude this section, we discuss the connection between these
results and Theorem B. Let $X$ be a special $\Pi$-manifold
over $M$, and let $\cale (X, \calf)$ denote the set of
bundle isomorphism classes of principal $(\Pi, G)$-bundles over $X$
with isotropy group system 
fine equivalent to
$\calf:=(H, U_A; \rho, \rho_A)$. Here a bundle
isomorphism is a $\Pi$-equivariant $G$-bundle isomorphism
$\phi\colon E_1 \to E_2$ covering the {\it identity\/} on $X$.
\begin{Lemma}
There is a commutative diagram: 
$$\xymatrix@C-2pt{\cale(X,\calf)\ar[r]\ar[d]^v
&\cals [H, U_A;\rho,\rho_A] 
\ar[r]^{\lambda}\ar[d]^{\approx}  &
\cals[H, U_A]\ar[d]^{\approx}
\\
\big [M, \partial M_A;BZ_\rho, BZ_{\rho_A}\big ]
\ar[r]^(0.45)u 
&\big [M, \partial M_A;B\Gamma\langle \rho\rangle, B\Omega_A\langle \rho_A \rangle\big ]
\ar[r] &
\big [M, \partial M_A;B\Gamma, B\Omega_A\big ]
}$$
where the horizontal composites
have images represented by the element $[X]$.
\end{Lemma}
\begin{proof}
The exact sequences
$1 \to Z_\rho\to\Gamma\langle\rho\rangle\to \Gamma$
and
$1 \to Z_{\rho_\alpha} \to\Omega_\alpha\langle\rho\rangle
\to \Omega_\alpha$
 of groups
induce a map 
$u\colon \big [M, \partial M_A; BZ_\rho,BZ_{\rho_A}\big ]
\to 
\big [M, \partial M_A;B\Gamma\langle \rho\rangle, B\Omega_A\langle \rho_A \rangle\big ]$.
By Theorem \ref{bundles} and Theorem \ref{special},
we also have a {\it surjective\/} map
$$v\colon \cale(X, \calf) \to 
\big [M, \partial M_A; BZ_\rho,BZ_{\rho_A}\big ]$$
induced by $u$ and our construction of equivariant
bundles.
\end{proof}
\begin{Theorem}\label{bundlesX}
Let $X$ be a special $\Pi$-manifold over $M$,
and $\calf=(H, U_A;\rho, \rho_A)$ an admissible
isotropy group system. Then
$$v\colon \cale(X, \calf) \cong 
\big [M, \partial M_A; BZ_\rho,BZ_{\rho_A}\big ]\ .$$
\end{Theorem}
\begin{proof}
The map $v$ is given in the diagram above, and we
have already observed that it is surjective.
Suppose that $\xi_1$, $\xi_2 \in \cale(X, \calf)$
with $v(\xi_1) = v(\xi_2)$.
Then we have a continuous map
$$(M\times I, \partial (M\times I))\to 
(BZ_\rho,BZ_{\rho_A})$$
realizing the homotopy between the classifying
maps for $\xi_1$ and $\xi_2$.
By the surjectivity of $v$ for $(\Pi, G)$-bundles
over $M\times I$, we get a bundle $(E, p)$
over $X\times I$ which restricts to $\xi_1$
and $\xi_2$ at the ends $X\times \partial I$.
Since $E \cong E_0\times I$, we get
 $\xi_1 \cong \xi_2$.
\end{proof}

Let $\aut (X)$ be the group of $\Pi$-equivariant isotopy classes
of $\Pi$-equivariant diffeomorphisms
of $X$ over the identity of $M$. The group
$\aut(X)$ preserving
 acts on 
$\cale (X, \calf)$
 by pulling back: $f\cdot \xi:=(f^\mun)^*\xi$. 
This action is well defined by
the equivariant Covering Homotopy Theorem of Palais \cite[II.7.3]{B}. 
If 
$$\lambda\colon \cals [H, U_A; \rho, \rho_A] \to
\cals [H, U_A]$$
denote the natural forgetful map, then there is an induced map
$$\psi\colon\cale (X, \calf) \to \lambda^\mun(X)$$
given by applying our stronger equivalence relation on bundles
(which allows $\phi$ to cover a self-diffeomorphism of $X$).

\begin{Proposition}\label{bundlesXaut}
The map $\psi$ induces  a bijection between $\lambda^\mun(X)$ and 
the quotient of
 $\cale (X, \calf)$ by the action of $\aut(X)$.
\end{Proposition}
\begin{proof}
The map  from one set of bundles to the other
is defined by regarding a $(\Pi, G)$-bundle over
$X$ as an element of $\lambda^\mun(X)$, and
this is well-defined since
 the equivalent relation in $\cals [H, U_A; \rho, \rho_A]$
 is stronger.
Moreover, two bundles with isotropy group system
$\calf$ over $X$ are equivalent if and only
if they are in the same orbit of the action of
 $\aut(X)$,
hence our correspondence is injective.
On the other hand, if $(E',p')$ is a bundle with base space
$X'$ in  $\lambda^\mun(X)$,
then there exists a $\Pi$-equivariant diffeomorphism
$h\colon X \to X'$ covering the identity on $M$.
Then $E:=h^*(E')$ is an equivalent element in $\lambda^\mun(X)$,
and is a bundle over $X$,
so our correspondence is surjective.
\end{proof}

These results and Theorem \ref{bundles} can sometimes
be used for explicit
classification of equivariant $(\Pi, G)$-bundles over special 
$\Pi$-manifolds. Notice that Bredon in \cite[V.7]{B} together with
\cite[Theorem V.6.4]{B}
has determined $\aut (X)$ in many cases of interest.

\medskip
 We shall now specialize to special 
 $\Pi$-manifolds over $I:=[-1,1]$, and extend Theorem B to this setting.
Examples include cohomogeneity 1 actions on spheres, classified
in \cite[Thm.~C, Table~II]{Straume}.
 The two components $\{\pm 1\}$ of
$\partial I$ are denoted $\{\pm\}$, and the notation $\calf=(H,U_\smallpm;\rho,\rho_\smallpm)$ will be used
for the admissible $(\Pi,G)$-isotropy group systems, as well as
 $\Gamma$, $\Omega_\smallpm$, etc. The classification of 
special $\Pi$-manifolds over $I$ takes the following form.

\begin{Theorem}\label{Ispecial} 
Let $(H, U_\smallpm)$ be an admissible isotropy group system over $I$.
Then ${\mathcal S}[H, U_\smallpm]$ is in bijection with the  double cosets
$\pi_0(\Omega_\smallminus)\backslash\pi_0(\Gamma)/\pi_0(\Omega_\smallplus )$.
\end{Theorem}

\begin{proof} This follows directly from Theorem \ref{special}, since
$$
[I,\partial I;B\Gamma,B\Omega_\smallpm]
\cong 
\pi_1(B\Omega_\smallminus)\backslash\pi_1(B\Gamma)/\pi_1(B\Omega_\smallplus)
\cong
\pi_0(\Omega_\smallminus)
\backslash\pi_0(\Gamma)/\pi_0(\Omega_\smallplus ).
\hfill \qed$$
\renewcommand{\qed}{}
\end{proof}
\begin{Remark}\label{rkIspecial}
The bijection of Theorem \ref{Ispecial} can be seen in a more
constructive way. 
Given a special $\Pi$-manifold $X$ over $I$, we can
choose an {\it $H$-meridian} $c\colon  I\to X$, i.e. a smooth section 
of $X\to I$ so that $\Pi_{c(t)}=H$ for $t$ in the interior of $I$; 
this can be obtained from a smooth section of 
the (trivial) principal $\Gamma$-bundle $P\to I$. Let
$U_\smallpm:=\Pi_{ c(\pm 1)}$.
We say that $(H,U_\smallpm)$ is an 
{\it $H$-meridian isotropy group system} for $X$.
Choosing another smooth section of $P$ gives isotropy groups
 conjugate
to $U_\smallpm$ by elements in the same connected component of $N(H)$.

As in (\ref{xialph}), the special $\Pi$-manifold $X$
can be reconstructed as a quotient of $I\times\Pi$:
$$ X= (I\times \Pi/H)\,\big / \{(\pm 1,g)\sim (\pm 1,gu_\pm),\,
\forall\, u_\smallpm\in  U_\smallpm\}$$
Therefore, $X$ is determined by the subgroups $U_\smallpm$
(compare \cite[Prop. 1.6]{GZ}). Moreover, any set $(H, U'_\smallpm)$,
where the $U'_\smallpm$ are conjugate to $U_\smallpm$,
occurs as an $H$-meridian isotropy group system for some special
$\Pi$-manifold $X'$ over $I$ (proved as in (\ref{pfCsur})).
In this way,
${\mathcal S}[H, U_\smallpm]$ is a quotient of $\pi_0(N(H))\times\pi_0(N(H))$.
The diagonal group acts trivially,  and by carefully examining 
the relevant equivalence relation
one sees that ${\mathcal S}[H, U_\smallpm]$ is in bijection with 
$\pi_0(\Omega_\smallminus )\backslash\pi_0(\Gamma)/\pi_0(\Omega_\smallplus )$. 
\end{Remark}

Let $X$ be a special $\Pi$-manifold over $I$.
We will now compute
the set $\cale(X,G)$ of isomorphism classes of $\Pi$-equivariant
principal $G$-bundles over $X$.
Fix an $H$-meridian isotropy group system 
$(H, U_\smallpm)$ for $X$
induced by a smooth $H$-meridian $c\colon I \to X$.
If $\xi = (E, p)$ is a $(\Pi,G)$-bundle over $X$, 
then we can choose
an isotropic lift $\tilde c\colon I \to E$, and obtain
an admissible isotropy group
system $\calf =(H,U_\smallpm;\rho,\rho_\smallpm)$ 
 for the bundle.

 Since the composition of an isotropic lift with a bundle isomorphism
is again isotropic, 
the conjugacy classes $[\rho_\smallpm]\in \calr(U_\smallpm,G)$ depend only on
$[\xi]\in \cale(X,G)$. 
This  defines a map 
$$J\colon  \cale(X,G)\to
\calr(U_\smallminus,G)\times
\calr(U_\smallplus ,G)\ .$$
We write $\calr(U_\smallminus,G)\times_H\calr(U_\smallplus ,G)$ for the set 
of pairs $([\rho_\smallminus ],[\rho_\smallplus ])
\in\calr(U_\smallminus ,G)\times\calr(U_\smallplus ,G)$
such that $\res [\rho_\smallminus ] = \res [\rho_\smallplus ]$ in $\calr(H,G)$.

Theorem B generalizes to special manifolds over $I$
as follows:

\begin{Theorem}\label{thb'}
Let $X$ be a special $\Pi$-manifold over $I$ realizing the isotropy
group system $(H, U_\smallpm)$.
The set $\cale(X,G)$ of isomorphism classes of $(\Pi, G)$-bundles
over $X$ is determined by the following
properties:
\begin{enumerate}
\renewcommand{\labelenumi}{(\roman{enumi})}
\item the image of $J$ is 
$\calr(U_\smallminus ,G)\times_H\calr(U_\smallplus ,G)$.
\item Let $\rho_\smallpm \colon U_\smallpm\to G$ be two smooth homomorphisms
such that $\rho_\smallminus \, {\scriptstyle |}_H=
\rho_\smallplus \, {\scriptstyle |}_H=\colon \rho$. Then there is a bijection
between
$J^\mun([\rho_\smallminus ],[\rho_\smallplus ])$ and
 the set of double cosets 
$\pi_0(Z_{\rho_\smallminus })\backslash\pi_0(Z_\rho)/\pi_0(Z_{\rho_\smallplus })$.
\end{enumerate}
\end{Theorem}

\begin{proof} Since $\rho_\smallpm$ comes from an admissible system,
the image of $J$ is contained in $\calr(U_\smallminus ,G)\times_H\calr(U_\smallplus ,G)$.
On the other hand, let $\rho_\smallpm \colon U_\smallpm\to G$ be two smooth homomorphisms
such that $\rho_\smallminus \, {\scriptstyle |}_H=
\rho_\smallplus \, {\scriptstyle |}_H=\colon \rho$.
The special $(\Pi\times G)$-manifold constructed as in 
Remark \ref{rkIspecial}, with isotropy group system 
$(H\langle\rho\rangle,U_\smallpm\langle\rho_\smallpm\rangle)$ 
(associated to an $H\langle\rho\rangle$-meridian) is a $(\Pi,G)$-bundle $\xi$
over $X$ (the special $\Pi$-manifold with $H$-meridian isotropy 
group system
$(H,U_\smallpm)$). One has 
$J(\xi)=([\rho_\smallminus ],[\rho_\smallplus ])$, which 
proves Part (i). 

If $(E,p) \in 
J^\mun([\rho_\smallminus ],[\rho_\smallplus ])$
then its isotropy group system is fine
equivalent to $\calf = (H, U_\smallpm; \rho,
\rho_\smallpm)$. In the notation introduced
earlier, we have
$$\cale(X, \calf) = J^\mun([\rho_\smallminus ],[\rho_\smallplus ])$$
and the result now follows from Theorem \ref{bundlesX}.
\end{proof}

We also have a more explicit version of
Proposition \ref{bundlesXaut}.
Note that  $J(f\cdot\xi)=J(\xi)$, for $f \in 
\aut(X)$, so we must
 investigate the action of $\aut (X)$
on a pre-image 
$\cale(X, \calf) \cong J^\mun([\rho_\smallminus ],[\rho_\smallplus ])$.
The group $\aut(X)$ has a homotopy description:
choose a base point $(\bullet) \in \Omega_\smallminus \backslash\Gamma /\Omega_\smallplus $ which corresponds
to the class of $X$ in 
$\pi_0(\Omega_\smallminus )\backslash\pi_0(\Gamma) /\pi_0(\Omega_\smallplus )$.

\begin{Proposition}\label{autX}
{\rm (Bredon \cite[V.7.3, VI.6.4]{B})}
There is a group anti-isomorphism 
$$\aut(X)\cong [I,\partial I;
\Gamma,\Omega_\smallpm  ]_{\bullet}$$
\end{Proposition}

\preu The pointed maps $(I, \partial I) \to (\Gamma, \Omega_\smallpm)$ send $\partial I$ into the
component of $\Omega_\smallminus \backslash\Gamma /\Omega_\smallplus$ containing the base point $(\bullet)$.
Let $f\in\aut(X)$. Using our $H$-meridian $c$, one defines a smooth
path $d\colon I\to \Pi/H$ by the formula
$$ f(c(t)) = d(t)\cdot c(t).$$
The map $f$ being $\Pi$-equivariant, one has, for all $h\in H$
$$hd(t)\cdot c(t)=h\cdot f(c(t))=f(hc(t))=f(c(t)) = d(t) c(t)$$
This implies, for all $t\in I$, that $d(t)\in  N(H)/H=\Gamma$.
For $t=\pm 1$, one gets in addition that $d(\pm 1)\in\Omega_\smallpm$.
We check that this defines an anti-homomorphism from $\aut(X)$ to 
the group $[I,\partial_\smallminus I,\partial_\smallplus I;\Gamma,\Omega_\smallminus ,\Omega_\smallplus ]$.
Now, if $d(t)$ represents a class in the latter, the formula 
$$f_d(\alpha c(t)):=\alpha d(t) \cdot c(t) \quad ,\quad \alpha\in\Pi, t\in I$$
defines an element $f_d$ of $\aut(X)$ and constitutes an inverse to
the above anti-homomorphism. \cqfd

Suppose that the homomorphisms $\Gamma\langle\rho\rangle\to\Gamma$
and $\Omega_\smallpm\langle\rho\rangle\to\Omega_\smallpm$ are surjective.
Then
we have a fibre bundle
\eqncount
\begin{equation}\label{fibZZ}
Z_{\rho_\smallminus }\backslash
 Z_\rho/Z_{\rho_\smallplus } \to
\Omega_\smallminus \langle\rho\rangle\backslash
\Gamma\langle\rho\rangle /
\Omega_\smallplus \langle\rho\rangle \to 
\Omega_\smallminus \backslash
\Gamma /
\Omega_\smallplus  .
\end{equation}
Therefore, $\pi_1(\Omega_\smallminus \backslash\Gamma /\Omega_\smallplus )$ acts 
on $\pi_0(Z_{\rho_\smallminus }\backslash Z_\rho/Z_{\rho_\smallplus })
=\pi_0(Z_{\rho_\smallminus })\backslash \pi_0(Z_\rho)/\pi_0(Z_{\rho_\smallplus })$.

\begin{Theorem}\label{thb'quotients}
Suppose that the homomorphisms $\Gamma\langle\rho\rangle\to\Gamma$
and $\Omega_\smallpm\langle\rho\rangle\to\Omega_\smallpm$ are surjective.
Let $\rho_\smallpm \colon U_\pm\to G$ be two smooth homomorphisms
such that $\rho_\smallminus \, {\scriptstyle |}_H=
\rho_\smallplus \, {\scriptstyle |}_H=\colon \rho$. Then the quotient of
$J^\mun([\rho_\smallminus ],[\rho_\smallplus ])$ by the action of $\aut(X)$
is in bijection with with the quotient of 
$\pi_0(Z_{\rho_\smallminus })\backslash \pi_0(Z_\rho)/\pi_0(Z_{\rho_\smallplus })$
by the action of $\pi_1(\Omega_\smallminus \backslash\Gamma /\Omega_\smallplus )$.
\end{Theorem}

\preu  
By Theorem \ref{bundlesX},  there is a surjective
map 
$\psi\colon J^\mun([\rho_\smallminus ],[\rho_\smallplus ])\to \lambda^\mun([X])$.
Since $\Gamma\langle\rho\rangle\to\Gamma$
and $\Omega_\smallpm\langle\rho\rangle\to\Omega_\smallpm$ are surjective,
we have the fibration \equationref{fibZZ}. From the
homotopy exact sequence of this fibration, 
we see that $\lambda^\mun([X])$ is in bijection with 
the quotient of 
$\pi_0\big(Z_{\rho_\smallminus }\backslash Z_\rho/Z_{\rho_\smallplus }\big)$
by the action of 
$\pi_1(\Omega_\smallminus \backslash\Gamma /\Omega_\smallplus )$.
But this is just the action of $\aut(X)$ by
Proposition \ref{autX}.
\end{proof}

\begin{Example}\label{ex1}
 Consider the standard action of $\Pi=SO(n)$
on $S^n$, with $H=SO(n-1)$. As $U_\smallpm=\Pi$, the homomorphisms
$\Gamma\langle\rho\rangle\to\Gamma$
and $\Omega_\smallpm\langle\rho\rangle\to\Omega_\smallpm$ are surjective.

When $n=2$, 
$H$ is the trivial group. If $G$ is connected,
then $Z_\rho=G$ is connected and the map $J$ is a bijection,
as seen in Theorem A. Also, the group of $\Pi$-automorphism
of $S^2$ is equal to $\Pi$, so $\aut(S^2)$ is trivial.

If $n\geq 3$,
the group $\Gamma=\Omega_\smallpm$ has 2 elements,
the non-trivial one represented by the diagonal matrix $D:={\rm
Diag\,}(1,\dots 1,-1,-1)$.
The space $\Omega_\smallminus \backslash\Gamma/\Omega_\smallplus $ being reduced to a point,
it follows from Theorem \ref{thb'quotients} that the group $\aut(S^n)$
acts trivially on $\cale(S^n,G)$. This has the following consequence:

\begin{Proposition}
Let $\rho_\smallpm\colon SO(n)\to G$ be two representations into a compact Lie group $G$
such that $\rho_\smallpm{}_{|SO(n-1)}=\rho$. Then, the element
$\rho_\smallminus (D)\rho_\smallplus (D)^\mun$ belongs to 
$Z_\rho$ and represent the trivial 
element in 
 $\pi_0(Z_{\rho_\smallminus })\backslash
\pi_0(Z_\rho)/\pi_0(Z_{\rho_\smallplus })$.
\end{Proposition}

\preu
As $\rho_\smallpm{}_{|H}=\rho$, one has, for all $h\in H$,
$$\rho_\smallplus (D)^\mun\rho(h)\rho_\smallplus (D)=\rho_\smallminus (D)^\mun\rho(h)\rho_\smallminus (D).$$
Therefore, $\rho_\smallminus (D)\rho_\smallplus (D)^\mun\in Z_\rho$.

Recall that our $H$-meridian for $S^n$ is 
$c(t)=(0,\dots,\cos(\pi t/2),\sin(\pi t/2))$. 
By Proposition \ref{TCnew} and its proof there exists 
a $(\Pi,G)$-bundle $\xi$ over $S^n$, with a 
$\rho$-isotropic lifting $\tilde c \colon  I\to E(\xi)$
of $c$ such that the isotropy representations associated
to $\tilde c(\pm 1)$ are $\rho_\pm$.
The curve $\tilde c$ is the horizontal 
lifting of $c$ for a $\Pi$-invariant connection. 
There is no problem to extend the definitions of
$c(t)$ and $\tilde c(t)$ for $t\in\bbr$. This defines $\mu\in G$ by
$$ \tilde c(3) = \tilde c(-1)\cdot \mu \ \hbox{ and }\ 
\tilde c(5) = \tilde c(1)\cdot \mu .$$
For $t\in [-1,1]$, one has $D\cdot\tilde c(t) = \tilde c(2-t)\rho_\smallplus (D)$
(since this is true for $t=1$ and both side are
horizontal). For $t=-1$, this gives
$$\tilde c(-1)\mu\rho_\smallplus (D)=\tilde c(3)\rho_\smallplus (D)
=D\cdot\tilde c(-1)=\tilde c(-1)\rho_\smallminus (D)$$
Therefore, $\mu=\rho_\smallminus (D)\rho_\smallplus (D)^\mun$.

Now, let $\xi_1:=\delta^*\xi$; one has 
$$E(\xi_1)=\{(x,u)\in S^n\times E(\xi)\mid \delta(x)=\pi(u)\}$$
and the $\Pi$-action on $E(\xi_1)$ comes from the diagonal action. 
A $\rho$-isotropic lifting $\tilde c_1 \colon  I\to E(\xi_1)$ of $c$ 
is then given by 
$$\tilde c_1(t):= (c(t),\tilde c(t-2)).$$
By Proposition \ref{TCnew}, one has $\xi=\xi_1$ in
$\cale(S^n,G)$ iff $\tilde J_\rho(\xi)=\tilde J_\rho(\xi_1)$ in
$\calr_\rho(n,G)$. By \ref{deftiga}, using the curves 
$\tilde c$ and $\tilde c_1$, one has in $\calr_\rho(n,G)$
$$\tilde J_\rho(\xi)=[\rho_\smallminus ,\rho_\smallplus ] \ \hbox{ and }\ 
\tilde J_\rho(\xi_1)=[\mu^\mun\rho_\smallminus \mu,\rho_\smallplus ].$$
This proves that $\mu$ represents the unit element in 
$\pi_0(Z_{\rho_\smallminus })\backslash \pi_0(Z_{\rho})/\pi_0(Z_{\rho_\smallplus })$. \cqfd
\end{Example}

\begin{Example}\label{ex2}
Consider the action of $SO(3)$ on the 
traceless $(3\times 3)$-symmetric matrixes by conjugation. Restricting 
this action to the unit sphere 
(for the invariant scalar product ${\rm tr\,}(AB)$)
makes $S^4$ a special $SO(3)$-manifold over $I$ 
(see \cite[\S 3]{GZ}) for more details on this classical cohomogeneity one
$SO(3)$-action on $S^4$). Here $H=S(O(1)\times O(1)\times O(1))$,
$U_\smallminus =S(O(1)\times O(2))$ and $U_\smallplus =S(O(2)\times O(1))$. Observe that
any homomorphism $\rho_\smallpm\colon U_\smallpm\to G$ is trivial unless its restriction
to $H$ is injective.

Take $G=SO(3)$ and use the standard inclusions of $U_\smallpm$
into $SO(3)$. Any non-trivial homomorphism $\rho_\smallminus $ is conjugate to 
$\tilde\rho_\smallminus \colon (\varepsilon,r_\theta)\mapsto (\varepsilon,r_\theta^p)$ and 
any non-trivial $\rho_\smallplus $ is conjugate to 
$\tilde\rho_\smallplus \colon (r_\theta,\varepsilon)\mapsto (r_\theta^q,\varepsilon)$,
where $p,q\in\bbn_{\rm odd}$, the set of positive odd integers.
Both $\tilde\rho_\smallpm$ restrict on $H$ to the identification 
$\rho$ of $H$ with the diagonal matrices of $SO(3)$. 
One has $Z_{\rho}=H$, 
$Z_{\tilde\rho_\smallminus} =\{{\rm diag\,}(1,1,1),{\rm diag\,}
(1,-1,-1)\}$ and
$Z_{\tilde\rho_\smallplus} =\{{\rm diag\,}(1,1,1),{\rm diag\,}
(-1,-1,1)\}$.  Therefore,
the set of double
cosets
$\pi_0(Z_{\tilde\rho_\smallminus })
\backslash\pi_0(Z_\rho)/\pi_0(Z_{\tilde\rho_\smallplus })$
is trivial 
and $\cale([S^4\to I],SO(3))$ is in bijection, via the map $J$ 
of theorem \ref{thb'}, with 
$\{0,0\}\cup\,\bbn_{\rm odd}\times \bbn_{\rm odd}$. 

In \cite[\S\ 3]{GZ}, $(SO(3),SO(3))$-bundles $P_{p,q}$ over $S^4$ 
are constructed for $p,q\in\bbz$ with $p\equiv 3 ({\rm mod\,} 4)$ 
(those come from $(S^3\times S^3)$-bundles). These bundles satisfy
$J(P_{p,q})=(|p|,|q|)$, so, up to isomorphism, the sign of $p$ and $q$
does not matter.
\end{Example}
\begin{Example}\label{ex3}
This is the complex analogue of Example \ref{ex2}.
One considers the action of $SU(3)$ on the 
traceless $(3\times 3)$-Hermitian matrixes by conjugation 
and restrict it
to the unit sphere. One thus gets a special $SU(3)$-manifold over $I$ 
diffeomorphic to $S^7$. The isotropy groups are 
$H=S(U(1)\times U(1)\times U(1))$,
$U_\smallminus=S(U(1)\times U(2))$ and $U_\smallplus =S(U(2)\times U(1))$. As in Example 1,
any homomorphism $\rho_\smallpm\colon U_\smallpm\to G$ is trivial unless its restriction
to $H$ is injective.

For $G=SU(3)$, one checks that any non trivial homomorphism of
$U_\smallpm$ to $G$ is conjugate to the standard inclusion. As $Z(H)=H$
in $SU(3)$, the map $J$ is injective and $\cale(S^7,SU(3))$
consists of two elements.
\end{Example}

\bibliographystyle{amsplain}
\bibliography{hh}
\end{document}